# A general and fast convolution-based method for peridynamics: applications to elasticity and brittle fracture

Siavash Jafarzadeh[1], Farzaneh Mousavi[1], Adam Larios[2,*], Florin Bobaru[1,*]

[1] *Department of Mechanical and Materials Engineering, University of Nebraska-Lincoln, Lincoln, NE 68588-0526, USA*
[2] *Department of Mathematics, University of Nebraska-Lincoln, Lincoln, NE 68588-0130, USA*

**Abstract**

We introduce a general and fast convolution-based method (FCBM) for peridynamics (PD). Expressing the PD integrals in terms of convolutions and computing them by fast Fourier transform (FFT), we reduce the computational complexity of PD models from O($N^2$) to O($N\log_2 N$), with $N$ being the total number of discretization nodes. Initial neighbor identification and storing neighbor information is not required, and, as a consequence, memory allocation scales with O($N$) instead of O($N^2$), common for existing methods. The method is applicable to bounded domains with arbitrary shapes and boundary conditions via an "embedded constraint" (EC) approach. We explain the FCBM-EC formulation for certain bond-based and state-based, linear and nonlinear PD models of elasticity and dynamic brittle fracture, as applications. We solve a 3D elastostatic problem and show that the FCBM reduces the computational time from days to hours and from years to days, compared with the original meshfree discretization for PD models. Large-scale computations of PD models are feasible with the new method, and we demonstrate its versatility by simulating, with ease, the difficult problem of multiple crack branching in a brittle plate.

1. Introduction

Peridynamics (PD) is a nonlocal extension of classical continuum mechanics used with great effect in modelling of fracture and damage [1, 2]. In contrast with classical (local) models described by partial differential equations (PDEs), governing equations in PD are integro-differential equations. In PD, nonlocal interactions are considered between points within a certain distance, replacing the usual spatial derivatives with integral operators. Because of this, the usual requirements of continuity and smoothness for displacements are eliminated, allowing one to model the emergence and evolution of discontinuities, such as cracks and damage, as natural parts of the solution to the governing equations [3, 4].

However, when usual discretization methods are used, nonlocality increases the cost for a PD model, relative to the cost of numerically approximating local models. The meshfree method with one-point Gaussian quadrature [5] (from now on, in this study we refer to it as "meshfree PD") and finite element methods (FEM) [6-9] have been used to compute numerical solutions to PD equations. For such equations, for a fixed finite range of nonlocal interaction (the horizon size), these methods scale as ($N^2$), where $N$ is the total number of nodes in the discretization. For large scale problems in 3D, the cost is prohibitive, even on massively parallel computers.

Various attempts have been made to reduce the cost of peridynamic simulations. Coupling the local theory with PD is one approach that uses a local model for parts of the domain, and uses the PD model

---

only at locations near cracks/damage as necessary [10, 11]. This approach does not work well for problems in which damage/cracks are (or become) widely distributed throughout the domain, such as in problems like impact fragmentation, etc. [12, 13].

The natural *convolutional structure* of PD formulations can be exploited using Fourier transforms as recently shown in [14, 15]. Using Fast Fourier Transform (FFT) algorithms, the new scaling for PD computations drops to $O(N\log_2 N)$ [16, 17]. These methods, however, are restricted to periodic domains. Fourier spectral methods have been used for nonlocal Allen-Cahn equation [18], fractional-in space reaction diffusion [19], and peridynamic diffusion and wave operators [14, 20-22], all in periodic domains. Another class of $O(N\log_2 N)$ methods was also introduced for 1D and 2D problems in [23-25], but these are still restricted to simple geometries and certain horizon shapes. Additionally, extension of these methods to 3D is not straightforward [23-25].

In [26, 27], we have introduced a fast convolution-based method (FCBM) for peridynamic diffusion problems. Two different techniques have been proposed in order to adapt the approach to problem sets on arbitrary domains and with general boundary conditions: the *volume penalization (VP)* [26] and the *embedded constraint (EC)* [27] techniques. The method has been validated in 1D, 2D, and 3D against exact nonlocal (manufactured) solutions and against local FEM-based solutions for problems in bounded domains with Dirichlet and Neumann boundary conditions [26, 27]. A diffusion problem in a 3D domain with an insulated cutout with billions of degrees of freedom was solved, over tens of thousands of time steps, in a matter of hours on a single processor. The same problem, if solved using the meshfree PD method, would have required years [27].

In this study, we extend the FCBM to the PD formulation for elastic deformations and fracture. We describe the class of constitutive models that have convolutional structures and can benefit from the remarkable efficiency of the method. Moreover, we demonstrate how to setup nonlinear problems like dynamic fracture in a brittle material as a convolution structure amenable to the FCBM treatment. This work also aims to serve as guide to demonstrate how one can construct material models with convolutional structures, and how to discretize different types of PD governing equation with this new method. As numerical examples, we formulate the FCBM-EC for linearized bond-based and linearized state-based PD isotropic solids, a nonlinear rotation-invariant model example, and a dynamic fracture problem. We validate the proposed damage model by comparing our simulation results against the published fracture PD simulations that use a popular critical bond strain criterion [5] for bond failure and the meshfree PD method for discretization.

This article is organized as follows: first a brief review of peridynamics equations of motion is provided, as well as concepts and notations related to PD constitutive modeling. Next, we discuss the convolutional from of PD material models through several examples. We describe the Fourier-based method for solving PD problems in periodic domain first, and then its extension to general bounded domains with arbitrary volume constraints/boundary conditions, using the EC method. We demonstrate the efficiency and validity of the introduced method by a 3D static and a 2D dynamic brittle fracture problems.

2. **Peridynamics**

We briefly describe the peridynamic formulation for elasticity and brittle fracture. Peridynamics is a nonlocal extension of continuum mechanics that unifies the governing equations for continuous media, media with discontinuities, and discrete media [1]. In this theory, each spatial point can interact with other points within its neighborhood up to a finite size distance. The shape of this neighborhood is arbitrary. In applications, however, for simplicity, for a point denoted by its position vector $\boldsymbol{x}$, the finite size neighborhood is taken as a sphere in 3D (a disk in 2D, or a segment in 1D) centered at $\boldsymbol{x}$ of radius $\delta$ [2]. This neighborhood, denoted by $\mathcal{H}_x$, is referred to as the *horizon region of $\boldsymbol{x}$* and $\delta$ is called the *horizon size*. Spatial points located inside $\mathcal{H}_x$ are referred to as the *family* of $\boldsymbol{x}$ and are denoted by $\boldsymbol{x}'$. The *bond vector* for each family pair is defined as $= \boldsymbol{x}' - \boldsymbol{x}$. Fig. 1 shows the schematic of a peridynamic body and the horizon of a generic point in the body.

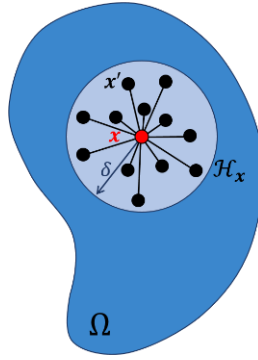

Fig. 1. A generic point $\boldsymbol{x}$, its horizon region $\mathcal{H}_x$, and its family points $\boldsymbol{x}'$ in a peridynamic body $\Omega$.

In this work, we employ the following notation to distinguish between different types of quantities: plain (non-bold) italic letters denote scalars, boldface italic letters denote vectors, and boldface non-italic letters denote tensors. Peridynamic states, discussed in section 2.2.2, are denoted by underlined letters.

### 2.1. Peridynamic Equation of motion

In continuum mechanics, the equation of motion can be expressed as:

$$\rho \frac{\partial^2 \boldsymbol{u}(\boldsymbol{x},t)}{\partial t^2} = \boldsymbol{L}(\boldsymbol{x},t) + \boldsymbol{b}(\boldsymbol{x},t) \tag{1}$$

where $\rho$ is mass density, $\boldsymbol{u}(\boldsymbol{x},t)$ is the displacement field at point $\boldsymbol{x}$ and time $t$, $\boldsymbol{L}(\boldsymbol{x},t)$ is the *internal force density*, and $\boldsymbol{b}$ is the *external body force density*. In the classical (local) theory, $\boldsymbol{L} = \boldsymbol{L}^\mathbf{C}(\boldsymbol{x},t) = \boldsymbol{\nabla} \cdot \boldsymbol{\sigma}$ (the divergence of the stress tensor) and the stress tensor dependency on strains/displacements is defined via a constitutive relationship. The superscript C denotes the classical definition of $\boldsymbol{L}$.

In peridynamics, $\boldsymbol{L}(\boldsymbol{x},t)$ is expressed with a nonlocal term instead: an integral which sums up all of the pair-wise forces between $\boldsymbol{x}$ and its family points $\boldsymbol{x}'$:

$$\boldsymbol{L}^{\mathbf{PD}}(\boldsymbol{x},t) = \int_{\mathcal{H}_x} \boldsymbol{f}(\boldsymbol{x},\boldsymbol{x}',t) \, dV_{x'} \tag{2}$$

The function $\boldsymbol{f}(\boldsymbol{x},\boldsymbol{x}',t)$ is called the *dual force density* and describes the net force between the unit volume at $\boldsymbol{x}$ and the unit volume at $\boldsymbol{x}'$ [1]. The specific form of $\boldsymbol{f}$ is determined by the chosen peridynamic constitutive model. The superscript PD denotes the peridynamics definition of $\boldsymbol{L}$. In the rest of the manuscript, we drop the superscript PD and use $\boldsymbol{L}$ to denote the PD internal force density.

As noted in the introduction, contrary to the local formulation, the PD operator $\boldsymbol{L}(\boldsymbol{x},t)$ has no continuity/smoothness requirements on the displacement field, allowing for cracks (as discontinuities in the displacement field) to naturally emerge and propagate in peridynamic solutions.

## 2.2. Peridynamic elastic and brittle fracture constitutive models in bond-based and state-based formulations

The original PD model was a "bond-based" model ([3]), which was later generalized to the state-based formulation of PD ([28]). Next, we review these formulations for describing elasticity and brittle failure.

### 2.2.1. Bond-based models

In bond-based materials, the dual force density between two family points depends *only* on deformation of their bond, and is independent of the deformation of other bonds in the same family. For example, in an elastic bond-based solid, dual force density depends only on the displacements of the family pair, i.e. $f = f(x, x', u, u', t)$ where $u' = u(x', t)$. Bond-based models suffer from certain restrictions on the material responses they produce. For example, Poisson ratio in bond-based isotropic elastic solids is restricted to a fixed number. To overcome such limitations, PD theory was extended to the state-based formulation by introducing the notion of *peridynamics states* [28].

### 2.2.2. State-based models

Peridynamics states are nonlinear (in general) operators, defined at a point and time, that map bonds into scalars or vectors, resulting in scalar-states or vector-states, respectively. A vector state is a generalization of the concept of a second-order tensor, which are linear mappings. Here, scalar-states are denoted by underlined italic lowercase letters, and vector-states are denoted by underlined boldface uppercase letters. The square bracket, $[\cdot]$, in front of a state encloses the point and the time at which the state is defined. The angle bracket $\langle \cdot \rangle$, if used, encloses the bond on which the state is operating. A state with the angle brackets in front of it then refers to the particular value associated with the enclosed bond. Here are some of the frequently used states in PD models:

$\underline{\mathbf{X}}[x]$, the *identity state* at point $x$ which returns back the bond vector for any given bond:

$$\underline{\mathbf{X}}[x]\langle \xi \rangle = \xi \tag{3}$$

$\underline{\mathbf{Y}}[x,t]$, the *deformation state* and returns the deformed bond vector at time $t$:

$$\underline{\mathbf{Y}}[x,t]\langle \xi \rangle = (x' + u') - (x + u) = \xi + \eta \tag{4}$$

where $\eta = u' - u$ is the relative displacement.

The *extension state* ($\underline{e}$) is a scalar state defining bond elongation:

$$\underline{e} = |\underline{\mathbf{Y}}| - |\underline{\mathbf{X}}|, \text{ i.e. } \underline{e}[x,t]\langle \xi \rangle = |\xi + \eta| - |\xi| \tag{5}$$

The *identity scalar-state* ($\underline{x}$) and the *deformation scalar-state* ($\underline{y}$) are used to defined the length of the bonds in reference and deformed configurations, respectively:

$$\underline{x} = |\underline{\mathbf{X}}|, \text{ and } \underline{y} = |\underline{\mathbf{Y}}|, \text{ i.e. } \underline{x}[x]\langle \xi \rangle = |\xi|, \text{ and } \underline{y}[x,t]\langle \xi \rangle = |\xi + \eta|. \tag{6}$$

$\underline{\mathbf{T}}[x,t]$ is the *force state* and returns the *bond force* $t$ for a given bond:

$$\underline{\mathbf{T}}[x,t]\langle \xi \rangle = t \tag{7}$$

where $t$ is the force vector that the unit volume at $x'$ exerts on a unit volume at $x$.

The dual force density $f$ in state-based PD is defined as follows:

$$f(x, x', t) = \underline{\mathbf{T}}[x, t]\langle\xi\rangle - \underline{\mathbf{T}}[x', t]\langle-\xi\rangle \tag{8}$$

According to Eq. (8), dual force density is equal to the force density that the unit volume at $x'$ exerts on $x$, minus the force density that unit volume at $x$ exerts on $x'$, i.e. it is the net force between the two unit volumes.

In classical continuum mechanics, constitutive models are tensorial equations, connecting a stress tensor to the strain tensor or the deformation gradient. In state-based PD, constitutive models similarly define the relationships between force states and deformation states. For example, an elastic state-based constitutive model is defined by $\underline{\mathbf{T}} = \underline{\mathbf{T}}(\underline{\mathbf{Y}})$, meaning that the force state is a function of the deformation state only. In state-based PD, the force density of a bond depends on the deformation of all other bonds inside the horizons of its two end-nodes [29]. Bond-based PD can be viewed as a special case of the state-based theory.

A frequent mathematical operation on PD states is the *inner product of the states*. Let $\underline{a}$ and $\underline{b}$ to denote two generic scalar states, and $\underline{\mathbf{A}}$ and $\underline{\mathbf{B}}$ to be two generic vector states. The inner product for scalar and vector states then are defined as follows [1]:

$$\underline{a} \bullet \underline{b} = \int_{\mathcal{H}_x} \underline{a}\langle\xi\rangle \underline{b}\langle\xi\rangle \, dV_{x'}, \qquad \underline{\mathbf{A}} \bullet \underline{\mathbf{B}} = \int_{\mathcal{H}_x} \underline{\mathbf{A}}\langle\xi\rangle \cdot \underline{\mathbf{B}}\langle\xi\rangle \, dV_{x'} \tag{9}$$

where ($\bullet$) denotes the inner product operation on PD states, and ($\cdot$) denotes the dot product of two vectors.

### 2.2.3. Modeling of damage and fracture in peridynamics

One common way to describe fracture and damage in a peridynamic model is using the "bond breaking" approach [5]. In this approach, a history dependent binary scalar function $\mu$ is inserted in the PD constitutive model:

$$\mathbf{L}(x, t) = \int_{\mathcal{H}_x} \mu(x, x', t) f(x, x', t) \, dV_{x'}, \tag{10}$$

where

$$\mu(x, x', t) = \begin{cases} 1 & (x' - x) \text{ bond is intact at time } t \\ 0 & (x' - x) \text{ bond broken at time } t \end{cases}. \tag{11}$$

Being "broken" or "intact" is determined by a certain failure criterion. The *critical bond strain* is one of such criterion [2, 5], which sets $\mu$ to 0 if the bond strain exceeds a critical value. Such critical values can be defined in several different ways [30]. In these models, cracks emerge as outcomes of cascading bond breaking events under loading. One way to monitor the evolution of damage in these models is the *damage index*, defined at a point:

$$d(x, t) = 1 - \frac{\int_{\mathcal{H}_x} \mu(x, x', t) \, dV_{x'}}{\int_{\mathcal{H}_x} dV_{x'}} \tag{12}$$

which counts the "number" of broken bonds relative to the total number of bonds for that point at a given time instant. Note that the damage index is not the "definition" of damage in PD, it is merely one way to represent damage in PD. Damage in these bond-breaking PD models is determined by the individual bond failure events, which is a much richer quantity than, for example, the scalar or even tensor variables defined in continuum damage mechanics (see [31]). In some sense, PD damage is a mapping, not necessarily continuous, from the vector space of bonds at a point $x$ to a vector space of dimension equal to the number of bonds at a node. In its discretized version, PD damage can be therefore considered a vector state. The damage index in Eq. (12) is a scalar representation that does not carry the complete information of PD damage. For example, consider these two cases: 1) a point has lost 50% of its bonds uniformly and symmetrically in all directions; and 2) a point has lost 50% of its bonds, but all in one side. In both cases $d = 0.5$, but in the first case the solid is uniformly degraded (distributed damage), and the second case is a split in solid (a through crack).

### 3. Obtaining convolutional structures for peridynamic models

Many linearized peridynamic models, featuring integral operators, can naturally be expressed via convolutions. Certain nonlinear models are also shown to have a convolutional structure [27, 32]. We observe that in general, the following form of a possibly nonlinear PD integrand (e.g. the dual force density in equation of motion):

$$f(\boldsymbol{x}, \boldsymbol{x}', t) = \sum_{n=1}^{p} a_n(\boldsymbol{x}, t) b_n(\boldsymbol{x}', t) c_n(\boldsymbol{x} - \boldsymbol{x}', t) \tag{13}$$

where $a_n, b_n$ and $c_n$ are functions, and $p$ is arbitrary positive integer, leads to a PD model that possesses a convolutional structure. Indeed, we can express:

$$\int_{\mathcal{H}_x} f(\boldsymbol{x}, \boldsymbol{x}', t) \, dV_{\boldsymbol{x}'} = \int_{\mathcal{H}_x} \left[ \sum_{n=1}^{p} a_n(\boldsymbol{x}, t) b_n(\boldsymbol{x}', t) c_n(\boldsymbol{x} - \boldsymbol{x}', t) \right] dV_{\boldsymbol{x}'} \tag{14}$$

$$= \sum_{n=1}^{p} a_n(\boldsymbol{x}, t) \int_{\mathcal{H}_x} c_n(\boldsymbol{x} - \boldsymbol{x}', t) b_n(\boldsymbol{x}', t) \, dV_{\boldsymbol{x}'} = \sum_{n=1}^{p} a_n(\boldsymbol{x}, t) \, [b_n * c_n](\boldsymbol{x}, t)$$

The operator $(*)$ denotes the convolution integral.

This observation allows for the development of efficient numerical techniques that utilize Fourier transforms. In this section we discuss setting up bond-based and state-based PD models in convolutional structures, for linear elastic materials. Moreover, we show how to obtain a convolutional structure for important nonlinear elastic models, and even for fracture problems. For the fracture case, we introduce a new damage model as an alternative to models that use the critical bond strain criterion for describing bond failure and does not have a convolutional structure.

#### 3.1. Linear elastic bond-based PD model

The force density for an elastic bond-based PD material can be expressed in the form [3]:

$$\boldsymbol{f}(\boldsymbol{x}, \boldsymbol{x}', t) = f(|\boldsymbol{\xi} + \boldsymbol{\eta}|, \boldsymbol{\xi}) \frac{\boldsymbol{\xi} + \boldsymbol{\eta}}{|\boldsymbol{\xi} + \boldsymbol{\eta}|} \tag{15}$$

Consider the well-known homogeneous isotropic elastic bond-based material [2] as an example:

$$\boldsymbol{f}(\boldsymbol{x}, \boldsymbol{x}', t) = \alpha \omega(|\boldsymbol{\xi}|)(|\boldsymbol{\xi} + \boldsymbol{\eta}| - |\boldsymbol{\xi}|) \frac{\boldsymbol{\xi} + \boldsymbol{\eta}}{|\boldsymbol{\xi} + \boldsymbol{\eta}|} \tag{16}$$

where $\alpha$ is a scalar and $\omega$ is a radially symmetric function of the bond, called the *influence function*, and is zero outside of the horizon: $\omega = 0$ for $|\boldsymbol{\xi}| > \delta$. Although this model is linear in terms of the bond extension, it is nonlinear in terms of displacements. The model in Eq. (16) can be linearized in terms of displacement if needed. According to the linearization carried out in [33], one can rewrite Eq.(16) for small displacements as:

$$f(x, x', t) = \alpha \omega(|\xi|) \frac{\xi \otimes \xi}{|\xi|^2} \cdot \eta = C(\xi) \cdot \eta \qquad (17)$$

where ($\otimes$) and ($\cdot$) denote respectively the tensor product and the inner product operators, and **C** is a tensor-valued symmetric function of the bond.

To obtain the convolutional form in Eq. (14), we opt to proceed with the indicial notation (including Einstein summation convention) for expressing vector and tensor quantities.

We substitute Eq. (17) into Eq.(2) to obtain the convolutional form for this model:

$$\begin{aligned} L_i(x, t) &= \int_{\mathcal{H}_x} f_i(x, x', t) \, dV_{x'} = \int_{\mathcal{H}_x} C_{ij}(x' - x)[u_j(x', t) - u_j(x, t)] dV_{x'} \\ &= \int_{\mathcal{H}_x} C_{ij}(x' - x) u_j(x', t) \, dV_{x'} - \left[\int_{\mathcal{H}_x} C_{ij}(x' - x) dV_{x'}\right] u_j(x, t) \\ &= \int_{\mathcal{H}_x} C_{ij}(x - x') \cdot u_j(x', t) \, dV_{x'} - P_{ij} u_j(x, t) \\ &= [C_{ij} * u_j](x, t) - P_{ij} u_j(x, t) \end{aligned} \qquad (18)$$

for $i, j = 1,2,3$ in 3D. Note that Eq. (18) is consistent with the general form in Eq. (14), because one can write $P_{ij} = C_{ij} * 1$. In this model, since we investigate homogeneous materials, **C** is only a function of $\xi$ and is independent of $x$. As a result **P** is a constant tensor. From a computational point of view, **P** can be computed once and stored as a preprocessing step. The linear term $P_{ij} u_j(x, t)$ is then just a multiplication at each node which is much cheaper than computing a volume integral. In case, due to the large problem size the one-time computation of **P** by direct quadrature becomes expensive, then one could instead use the convolutional form $P_{ij} = C_{ij} * 1$ to compute it by FFT, like we do for the other term in the convolutional form presented by Eq. (18).

This linearized model in Eq. (17) is based on the "small displacements" assumption, and therefore, cannot be used for large displacements, including rigid body rotations. In general, the nonlinear model in Eq. (16) is used more frequently [34-36] because it is valid for large rotations due to the geometrical nonlinearity of the model. We could not obtain the convolutional form for the nonlinear model in Eq. (16) directly. However, in the next section we introduce an alternative nonlinear bond-based model which has a convolutional structure, and is invariant under arbitrary rigid body rotations.

### 3.2. Nonlinear elastic bond-based PD model

Here, we introduce a particular nonlinear elastic bond-based model and show how to express it using multiple convolutions. Consider the bond-based constitutive model:

$$f(x, x', t) = \alpha \omega(|\xi|) |\xi| \left(s + \frac{3}{2} s^2 + \frac{1}{2} s^3\right) \frac{\xi + \eta}{|\xi + \eta|} \qquad (19)$$

where

$$s = \frac{|\xi + \eta| - |\xi|}{|\xi|} \qquad (20)$$

Is the *bond strain*. This nonlinear model is invariant to rigid body motions and approximates the model in Eq. (16) for $s \ll 1$. To express $L(x, t)$ in its convolutional form, we first simplify $f$ in terms of displacements and bond vectors:

$$f(x, x', t) = \alpha\omega(|\xi|)|\xi|\left(s + \frac{3}{2}s^2 + \frac{1}{2}s^3\right)\frac{\xi + \eta}{|\xi + \eta|} = \tag{21}$$

$$= \frac{\alpha\omega(|\xi|)|\xi|}{2}[(s+1)^2 - 1](s+1)\frac{\xi + \eta}{|\xi + \eta|}$$

$$= \frac{\alpha\omega(|\xi|)|\xi|}{2}\left[\left(\frac{|\xi + \eta|}{|\xi|}\right)^2 - 1\right]\frac{|\xi + \eta|}{|\xi|}\frac{\xi + \eta}{|\xi + \eta|}$$

$$= \frac{\alpha\omega(|\xi|)}{2|\xi|^2}(|\xi + \eta|^2 - |\xi|^2)(\xi + \eta) = \frac{\alpha\omega(|\xi|)}{2|\xi|^2}(2\xi \cdot \eta + \eta \cdot \eta)(\xi + \eta)$$

Again, we use the indicial notation for $f$ in Eq. (21):

$$f_i(x, x', t) = \frac{\alpha\omega(|\xi|)}{2|\xi|^2}(2\xi_j\eta_j + \eta_j\eta_j)(\xi_i + \eta_i) \tag{22}$$

$$= \frac{\alpha\omega(|\xi|)}{2|\xi|^2}(2\xi_j u'_j - 2\xi_j u_j + u'_j u'_j - 2u'_j u_j + u_j u_j)(\xi_i + u'_i - u_i)$$

$$= \frac{\alpha\omega(|\xi|)}{2|\xi|^2}(2\xi_i\xi_j u'_j - 2\xi_i\xi_j u_j + \xi_i u'_j u'_j - 2\xi_i u'_j u_j + \xi_i u_j u_j + 2\xi_j u'_j u'_i$$

$$- 2\xi_j u_j u'_i + u'_j u'_j u'_i - 2u'_j u_j u'_i + u_j u_j u'_i - 2\xi_j u'_j u_i + 2\xi_j u_j u_i$$

$$- u'_j u'_j u_i + 2u'_j u_j u_i - u_j u_j u_i)$$

We then can write $L(x, t)$ as:

$$L_i = \int_{\mathcal{H}_x} f_i \, dV_{x'} = \int_{\mathcal{H}_x} \frac{\alpha\omega(|x' - x|)}{2|x' - x|^2}[2(x_i - x'_i)(x_j - x'_j)u'_j - 2(x_i - x'_i)(x_j - x'_j)u_j \tag{23}$$

$$- (x_i - x'_i)u'_j u'_j + 2(x_i - x'_i)u'_j u_j - (x_i - x'_i)u_j u_j - 2(x_j - x'_j)u'_j u'_i$$

$$+ 2(x_j - x'_j)u_j u'_i + 2(x_j - x'_j)u'_j u_i - 2(x_j - x'_j)u_j u_i + u'_j u'_j u'_i$$

$$- 2u'_j u_j u'_i + u_j u_j u'_i - u'_j u'_j u_i + 2u'_j u_j u_i - u_j u_j u_i]dV_{x'}$$

with $i, j = 1,2,3$ in 3D. Let the tensor $\mathbf{C}$, the vector $\mathbf{a}$, and the scalar $c$ denote the following functions of $\xi$:

$$\mathbf{C}(\xi) = \frac{\alpha\omega(|\xi|)}{|\xi|^2}\xi \otimes \xi, \qquad \mathbf{a}(\xi) = \frac{\alpha\omega(|\xi|)}{|\xi|^2}\xi, \qquad c(\xi) = \frac{\alpha\omega(|\xi|)}{|\xi|^2} \tag{24}$$

and let $\mathbf{P}$, $\mathbf{q}$, and $p$ to be the integrals of $\mathbf{C}$, $\mathbf{a}$, and $c$ over the horizon respectively. Since $\mathbf{C}$, $\mathbf{a}$, and $c$ are functions of $\xi$ and independent of $x$, their integrals can be computed at any $x$. We can then compute the integrals once at $x = 0$ for example, where $\xi = x'$:

$$\mathbf{P} = \int_{\mathcal{H}_{x=0}} \mathbf{C}(x')dV_{x'}, \qquad \mathbf{q} = \int_{\mathcal{H}_{x=0}} \mathbf{a}(x')dV_{x'}, \qquad p = \int_{\mathcal{H}_{x=0}} c(x')dV_{x'} \tag{25}$$

Using the definitions in Eqs.(24) and (25), Eq. (23) can be reorganized as follows:

$$L_i = C_{ij} * u_j - \left(\int_{\mathcal{H}_x} C_{ij} dV_{x'}\right) u_j - a_i * (u_j u_j) + (a_i * u_j) u_j - \frac{1}{2}\left(\int_{\mathcal{H}_x} a_i dV_{x'}\right) u_j u_j - a_j \quad (26)$$

$$* (u_j u_i) + u_j(a_j * u_i) + (a_j * u_j) u_i - \left(\int_{\mathcal{H}_x} a_i dV_{x'}\right) u_i u_j + \frac{1}{2} c * (u_j u_j u_i)$$

$$- [c * (u_i u_j)] u_j + \frac{1}{2}(c * u_i) u_j u_j - \frac{1}{2}[c * (u_j u_j)] u_i + (c * u_j) u_j u_i$$

$$- \frac{1}{2} u_j u_j u_i \left(\int_{\mathcal{H}_x} c \, dV_{x'}\right)$$

$$= C_{ij} * u_j - P_{ij} u_j - a_i * (u_j u_j) + (a_i * u_j) u_j - \frac{1}{2} q_i u_j u_j - a_j * (u_j u_i)$$

$$+ u_j(a_j * u_i) + (a_j * u_j) u_i - q_i u_i u_j + \frac{1}{2} c * (u_j u_j u_i) - [c * (u_i u_j)] u_j$$

$$+ \frac{1}{2}(c * u_i) u_j u_j - \frac{1}{2}[c * (u_j u_j)] u_i + (c * u_j) u_j u_i - \frac{1}{2} u_j u_j u_i p$$

Eq. (26) shows the convolutional form of the **L** integral for the nonlinear model (rigid body rotation invariant) introduced in Eq. (19).

### 3.3. State-based linear elastic and isotropic PD solid

In this example we discuss a well-known state-based linear elastic model, known as the *linear isotropic peridynamic solid* [33]. The force state for this material is given by:

$$\underline{\mathbf{T}} = \left(\frac{3\kappa - 5\mu}{m}\right) \underline{\omega x} \vartheta \underline{\mathbf{M}} + \frac{15\mu}{m} \underline{\omega e} \, \underline{\mathbf{M}} \quad (27)$$

where $\kappa$ and $\mu$ are the *bulk* and the *shear moduli*, $\vartheta$ is the *nonlocal dilatation* defined at each point by:

$$\vartheta(\mathbf{x},t) = \frac{3}{m} \underline{\omega x} \bullet \underline{e} = \frac{3}{m} \int_{\mathcal{H}_x} \underline{\omega}\langle \xi \rangle \, \underline{x}\langle \xi \rangle \, \underline{e}\langle \xi \rangle \, dV_{x'} \, , \quad (28)$$

$m$ is a normalization factor:

$$m(\mathbf{x},t) = \underline{\omega x} \bullet \underline{x} = \int_{\mathcal{H}_x} \underline{\omega}\langle \xi \rangle \, \underline{x}\langle \xi \rangle \, \underline{x}\langle \xi \rangle \, dV_{x'} \, , \quad (29)$$

$\underline{\omega}$ is the influence function scalar state:

$$\underline{\omega}\langle \xi \rangle = \omega(|\xi|) \, , \quad (30)$$

**M** is a unit vector state, giving the direction of the deformed bond:

$$\underline{\mathbf{M}}\langle \xi \rangle = \frac{\underline{\mathbf{Y}}\langle \xi \rangle}{|\underline{\mathbf{Y}}\langle \xi \rangle|} = \frac{\xi + \eta}{|\xi + \eta|} \, , \quad (31)$$

and $\underline{e}$ is the extension state defined in Eq. (5). Note that we can express Eq.(27) in terms of $\xi$ and $\eta$ as follows:

$$\underline{\mathbf{T}}\langle \xi \rangle = \left(\frac{3\kappa - 5\mu}{m}\right) \omega(|\xi|) |\xi| \vartheta \frac{\xi + \eta}{|\xi + \eta|} + \frac{15\mu}{m} \omega(|\xi|)(|\xi + \eta| - |\xi|) \frac{\xi + \eta}{|\xi + \eta|}. \quad (32)$$

The linearized version of this material for small displacements is [33]:

$$\underline{\mathbf{T}}\langle \boldsymbol{\xi}\rangle = \left(\frac{3\kappa - 5\mu}{m}\right)\omega(|\boldsymbol{\xi}|)\vartheta\boldsymbol{\xi} + \frac{15\mu}{m}\omega(|\boldsymbol{\xi}|)\frac{\boldsymbol{\xi}\otimes\boldsymbol{\xi}}{|\boldsymbol{\xi}|^2}\cdot\boldsymbol{\eta} \qquad (33)$$

with

$$\vartheta = \frac{3}{m}\int_{\mathcal{H}_x}\omega(|\boldsymbol{\xi}|)\,\boldsymbol{\xi}\cdot\boldsymbol{\eta}\,\mathrm{d}V_{x'} \qquad (34)$$

and

$$m = \int_{\mathcal{H}_x}\omega(|\boldsymbol{\xi}|)\,|\boldsymbol{\xi}|^2\,\mathrm{d}V_{x'} \qquad (35)$$

Based on this constitutive model, the internal force density is:

$$\begin{aligned}
\boldsymbol{L}(\boldsymbol{x},t) &= \int_{\mathcal{H}_x}\boldsymbol{f}(\boldsymbol{x},\boldsymbol{x}',t)\,\mathrm{d}V_{x'} = \int_{\mathcal{H}_x}\{\underline{\mathbf{T}}[\boldsymbol{x},t]\langle\boldsymbol{\xi}\rangle - \underline{\mathbf{T}}[\boldsymbol{x}',t]\langle-\boldsymbol{\xi}\rangle\}\mathrm{d}V_{x'} \\
&= \int_{\mathcal{H}_x}\bigg\{\left(\frac{3\kappa-5\mu}{m}\right)\omega(|\boldsymbol{\xi}|)\vartheta\boldsymbol{\xi} + \frac{15\mu}{m}\omega(|\boldsymbol{\xi}|)\frac{\boldsymbol{\xi}\otimes\boldsymbol{\xi}}{|\boldsymbol{\xi}|^2}\cdot\boldsymbol{\eta} \\
&\quad - \left(\frac{3\kappa-5\mu}{m}\right)\omega(|-\boldsymbol{\xi}|)\vartheta'(-\boldsymbol{\xi}) + \frac{15\mu}{m}\omega(|-\boldsymbol{\xi}|)\frac{(-\boldsymbol{\xi})\otimes(-\boldsymbol{\xi})}{|-\boldsymbol{\xi}|^2}\cdot(-\boldsymbol{\eta})\bigg\}\mathrm{d}V_{x'} \\
&= \int_{\mathcal{H}_x}\bigg\{\left(\frac{3\kappa-5\mu}{m}\right)\omega(|\boldsymbol{\xi}|)(\vartheta+\vartheta')\boldsymbol{\xi} + 2\frac{15\mu}{m}\omega(|\boldsymbol{\xi}|)\frac{\boldsymbol{\xi}\otimes\boldsymbol{\xi}}{|\boldsymbol{\xi}|^2}\cdot\boldsymbol{\eta}\bigg\}\mathrm{d}V_{x'}
\end{aligned} \qquad (36)$$

where $\vartheta' = \vartheta(\boldsymbol{x}',t)$. Let

$$\mathbf{C}(\boldsymbol{\xi}) = \frac{30\mu}{m}\omega(|\boldsymbol{\xi}|)\frac{\boldsymbol{\xi}\otimes\boldsymbol{\xi}}{|\boldsymbol{\xi}|^2}, \qquad \boldsymbol{a}(\boldsymbol{\xi}) = \left(\frac{3\kappa-5\mu}{m}\right)\omega(|\boldsymbol{\xi}|)\boldsymbol{\xi}. \qquad (37)$$

Replacing $\boldsymbol{\xi}$ and $\boldsymbol{\eta}$ in Eq. (36) with $\boldsymbol{x}'-\boldsymbol{x}$ and $\boldsymbol{u}'-\boldsymbol{u}$ respectively, one gets:

$$\begin{aligned}
L_i &= \int_{\mathcal{H}_x}a_i(\boldsymbol{x}'-\boldsymbol{x})(\vartheta+\vartheta')\mathrm{d}V_{x'} + \int_{\mathcal{H}_x}C_{ij}(\boldsymbol{x}'-\boldsymbol{x})(u_j'-u_j)\mathrm{d}V_{x'} \\
&= \vartheta\int_{\mathcal{H}_x}a_i(\boldsymbol{x}'-\boldsymbol{x})\mathrm{d}V_{x'} - \int_{\mathcal{H}_x}a_i(\boldsymbol{x}-\boldsymbol{x}')\vartheta'\mathrm{d}V_{x'} + \int_{\mathcal{H}_x}C_{ij}(\boldsymbol{x}-\boldsymbol{x}')u_j'\mathrm{d}V_{x'} \\
&\quad + \left(\int_{\mathcal{H}_x}C_{ij}(\boldsymbol{x}'-\boldsymbol{x})\mathrm{d}V_{x'}\right)u_j \\
&= \vartheta\left(\int_{\mathcal{H}_x}a_i(\boldsymbol{x}'-\boldsymbol{x})\mathrm{d}V_{x'}\right) - (a_i * \vartheta) + C_{ij}*u_j \\
&\quad + \left(\int_{\mathcal{H}_x}C_{ij}(\boldsymbol{x}'-\boldsymbol{x})\mathrm{d}V_{x'}\right)u_j \\
&= \vartheta\left(\int a_i\mathrm{d}V_x\right) - (a_i*\vartheta) + C_{ij}*u_j + \left(\int C_{ij}\mathrm{d}V_x\right)u_j \\
&= \vartheta q_i - (a_i*\vartheta) + C_{ij}*u_j + P_{ij}u_j
\end{aligned} \qquad (38)$$

To use the numerical method presented in this work, all volume integrals involved in the PD model should have convolutional structure, otherwise the FFT cannot be used to compute them and the proposed method's efficiency is lost. Therefore, $\vartheta(\boldsymbol{x},t)$ needs to be written as a convolution as well:

$$\vartheta = \frac{3}{m}\int_{\mathcal{H}_x} \omega(|\xi|)\,\xi_i\eta_i\, dV_{x'} = \frac{-3}{m}\int_{\mathcal{H}_x} a_i(x-x')(u_i'-u_i)\, dV_{x'} \quad (39)$$
$$= \frac{-3}{m}\int_{\mathcal{H}_x} a_i(x-x')u_i'\, dV_{x'} - \frac{3}{m}\left(\int_{\mathcal{H}_x} a_i(x'-x)dV_{x'}\right)u_i$$
$$= \frac{-3}{m}[(a_i * u_i) + q_i u_i]$$

Eqs.(38) and (39), provide the convolutional form for the linearized PD state-based elastic model described in Eq. (33).

### 3.4. An energy-based failure model that leads to a convolution structure for fracture problems

The failure model in PD, based on bond-breaking once bonds reach a critical strain value, does have a convolution structure. To be able to use the fast convolution-based method in fracture problems, we introduce a new energy-based bond-failure model.

#### 3.4.1. Pointwise energy-based bond breaking

We propose the following bond brittle failure model based on the strain energy density at the bond's two end points:

$$\mu(x, x', t) = \begin{cases} 1, & \text{if } W(x,t) \text{ and } W(x',t) \leq W_c \\ 0, & \text{if } W(x,t) \text{ or } W(x',t) > W_c \end{cases} \quad (40)$$

where $W(x, t)$ denotes the strain energy density at point $x$ and time $t$ and $W_c$ is referred to as the *critical bond strain energy density*. The formula for $W$ depends on the constitutive model. For example, the strain energy density for the model in Eq.(16) is [33]:

$$W(x,t) = \frac{1}{2}\int_{\mathcal{H}_x} \alpha\omega(|\xi|)(|\xi+\eta|-|\xi|)^2\, dV_{x'} \quad (41)$$

and for the linearized model in Eq.(17) is [37]:

$$W(x,t) = \frac{1}{2}\int_{\mathcal{H}_x} \eta \cdot C(|\xi|) \cdot \eta\, dV_{x'} \quad (42)$$

Note that for models with bond-breaking rules, $\mu$ is inserted in the integrand of $W$ in formulas to leave out broken bonds.

Eq. (40) states that a bond is broken and no longer carries load if the strain energy density at either of its end points exceeds $W_c$. A pointwise interpretation of Eq. (40) is that once $W$ at a point reaches $W_c$, that point loses all of its bonds, and it is completely detached from the body.

To obtain the convolutional form for this brittle failure model, we define the *integrity index* $\lambda(x, t)$ at each point as:

$$\lambda(x,t) = \begin{cases} 1 & W(x,t) \leq W_c \\ 0 & W(x,t) > W_c \end{cases} \quad (43)$$

Using the integrity index, we can write:

$$\mu(x, x', t) = \lambda(x,t)\lambda(x',t) = \lambda\lambda' \quad (44)$$

Next we show how this new damage model leads to a convolutional structure.

### 3.4.2. Convolutional structure of the new damage model

In general, we show that a PD damage model with $L(x,t) = \int_{\mathcal{H}_x} \mu f \, dV_{x'}$ has a convolutional structure if $f$ and $\mu$ are both in the form described by Eq.(13):

$$L(x,t) = \int_{\mathcal{H}_x} \mu(x,x',t) f(x,x',t) \, dV_{x'} \tag{45}$$

$$= \int_{\mathcal{H}_x} \left[ \sum_{n=1}^{p} a_n(x,t) b_n(x',t) c_n(x-x',t) \right] \left[ \sum_{m=1}^{q} y_m(x,t) w_m(x',t) z_m(x-x',t) \right] dV_{x'} = \int_{\mathcal{H}_x} \left[ \sum_{n=1}^{p} \sum_{m=1}^{q} a_n y_m(x,t) b_n w_m(x',t) c_n z_m(x-x',t) \right] dV_{x'}$$

$$= \sum_{n=1}^{p} \sum_{m=1}^{q} a_n y_m(x,t) \int_{\mathcal{H}_x} b_n w_m(x',t) c_n z_m(x-x',t) \, dV_{x'}$$

$$= \sum_{n=1}^{p} \sum_{m=1}^{q} a_n y_m(x,t) [(b_n w_m) * (c_n z_m)](x,t)$$

The new $\mu = \lambda \lambda'$ is in the form of Eq. (13). As a result, it can be used with any of the described constitutive models in section 3. For the linearized bond-based model described in Section 3.1 for example, according to Eq. (10), Eq. (18) and Eq.(44) we can write:

$$L_i = \int_{\mathcal{H}_x} \lambda \lambda' f_i \, dV_{x'} = \lambda \int_{\mathcal{H}_x} \lambda' f_i \, dV_{x'} = \lambda \int_{\mathcal{H}_x} \lambda' C_{ij}(x-x')(u'_j - u_j) \, dV_{x'} \tag{46}$$

$$= \lambda \int_{\mathcal{H}_x} C_{ij}(x-x')\lambda' u'_j \, dV_{x'} - \lambda u_j \int_{\mathcal{H}_x} C_{ij}(x-x')\lambda' \, dV_{x'}$$

$$= \lambda [C_{ij} * (\lambda u_j)] - \lambda u_j [C_{ij} * \lambda]$$

The strain energy density for this material model can be written as:

$$W(x,t) = \frac{1}{2} \int_{\mathcal{H}_x} \lambda \lambda' \boldsymbol{\eta} \cdot \mathbf{C}(x'-x) \cdot \boldsymbol{\eta} \, dV_{x'} = \frac{1}{2} \int_{\mathcal{H}_x} \lambda \lambda' \eta_i C_{ij}(x'-x) \eta_j \, dV_{x'} \tag{47}$$

$$= \frac{1}{2} \int_{\mathcal{H}_x} \lambda \lambda' C_{ij}(x'-x)(u'_i - u_i)(u'_j - u_j) \, dV_{x'}$$

$$= \frac{1}{2} \int_{\mathcal{H}_x} \lambda \lambda' C_{ij}(x'-x)(u'_i u'_j - u'_i u_j - u_i u'_j + u_i u_j) \, dV_{x'}$$

$$= \frac{1}{2} \lambda \left[ \int_{\mathcal{H}_x} C_{ij}(x'-x)\lambda' u'_i u'_j \, dV_{x'} - u_j \int_{\mathcal{H}_x} C_{ij}(x'-x)\lambda' u'_i \, dV_{x'} \right.$$

$$\left. - u_i \int_{\mathcal{H}_x} C_{ij}(x'-x)\lambda' u'_j \, dV_{x'} + u_i u_j \int_{\mathcal{H}_x} C_{ij}(x'-x)\lambda' \, dV_{x'} \right]$$

$$= \frac{1}{2} \lambda \{ [C_{ij} * (\lambda u_i u_j)] - u_j [C_{ij} * (\lambda u_i)] - u_i [C_{ij} * (\lambda u_j)] + u_i u_j [C_{ij} * \lambda] \}$$

The damage index in Eq. (12) can also be written in terms of convolutions if needed, as shown below. By defining a *characteristic influence function* with its support equal to the horizon region:

$$\omega_0(\xi) = \begin{cases} 1 & |\xi| \leq \delta \\ 0 & |\xi| > \delta \end{cases}, \qquad (48)$$

one can write:

$$d(\boldsymbol{x},t) = 1 - \frac{\int_{H_x} \mu \mathrm{d}V_{x'}}{\int_{H_x} \mathrm{d}V_{x'}} = 1 - \frac{\int_{\mathbb{R}^n} \lambda \lambda' \omega_0(|\boldsymbol{x}'-\boldsymbol{x}|) \; \mathrm{d}V_{x'}}{\int_{\mathbb{R}^n} \omega_0(|\boldsymbol{x}'-\boldsymbol{x}|) \; \mathrm{d}V_{x'}} = 1 - \frac{\int_{\mathbb{R}^n} \lambda \lambda' \omega_0(|\boldsymbol{x}-\boldsymbol{x}'|) \; \mathrm{d}V_{x'}}{\beta_0} \qquad (49)$$

$$= 1 - \frac{\lambda \int_{-\infty}^{\infty} \lambda' \omega_0(|\boldsymbol{x}-\boldsymbol{x}'|) \mathrm{d}V_{x'}}{\beta_0} = 1 - \frac{\lambda(\lambda * \omega_0)}{\beta_0}$$

### 3.4.3. Calibrating the failure model to Griffith's critical energy release rate

To use this model for predicting fracture, the failure threshold $W_c$ in Eq. (40) needs to be calibrated to a measurable fracture property of the material like Griffith's critical fracture energy. Inspired by the calibration method for the critical bond strain criterion [5], we perform the calibration by considering a "through" planar crack in a body and find a relationship between the PD failure threshold (here $W_c$) and the critical energy release rate $G_0$, which is a measurable quantity.

A through crack is one that splits the body into two completely separated parts. For the pointwise criterion given by Eqs. (40) and (43), *a through crack contains a $\delta$-thick layer with points that lost all of their bonds*. In Fig. 2 we show schematically why this is the case.

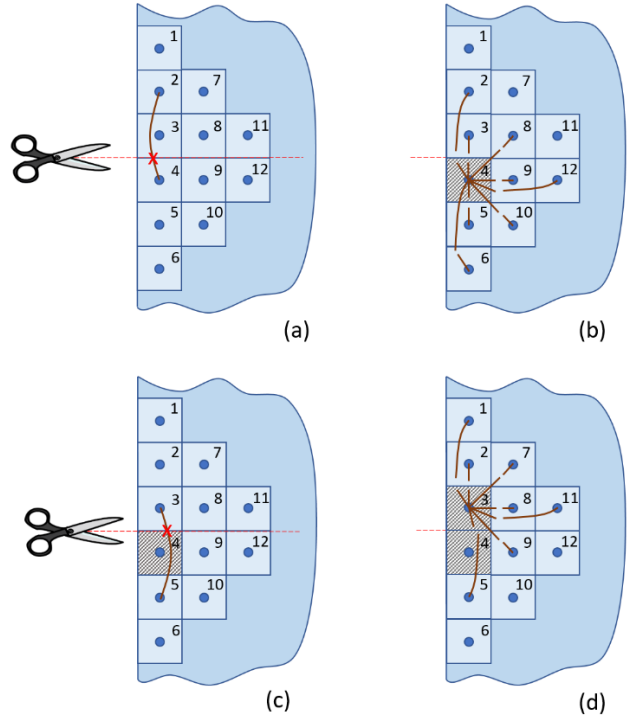

Fig. 2. Schematics of a through crack based on the new damage model in 2D, for a case with $\delta = 2\Delta x$. "Cutting" bonds along a line results in a $\delta$-thick layer of nodes with $d = 1$ (all of their bonds are broken). The figure shows the case for m = 2.

Fig. 2 shows what happens when a through crack is forming by breaking bonds along a line in 2D. To make the figure more readable, we show the case for $\delta = 2\Delta x$, which has fewer bonds to plot. In Fig. 2a, bond 2-4 is cut. According to Eq. (40) either node 2 or 4 has reached $W_c$ and therefore should lose *all* of its bonds. Fig. 2b shows the case for node 4. As shown, cutting the bond 2-4 resulted in breaking other bonds that are connected to node 4. As the through cutting along the line proceeds, we cut bond 3-5 (Fig. 2c). As a result either node 3 or 5 loses all of its bonds and reached the damage index of 1. Fig. 2d shows the case for node 3. We chose first node 4 and then node 3 for demonstration. All other possible scenarios in nodal failure resulting from cutting 2-4, 3-4, and 3-5 bonds that cross the through crack path leads to at least two adjacent failed nodes (d=1). This simple example shows how a through crack forms a $\delta$-thick layer of nodes with damage index of 1.

In general, a through crack with this model has a damage region with the thickness of $3\delta$: a $\delta$-thick layer in the middle with the damage index $d = 1$ and two $\delta$-thick layers on the sides with $0 < d < 0.5$. This is different from the case for the critical bond strain criterion where a through crack has a $2\delta$-thick damaged region with $0 < d < 0.5$ and the maximum $d$ being in the middle plane. Fig. 3 schematically compares a through crack in the new damage model with the one in the critical bond strain damage model in terms of thickness and damage profile across the crack.

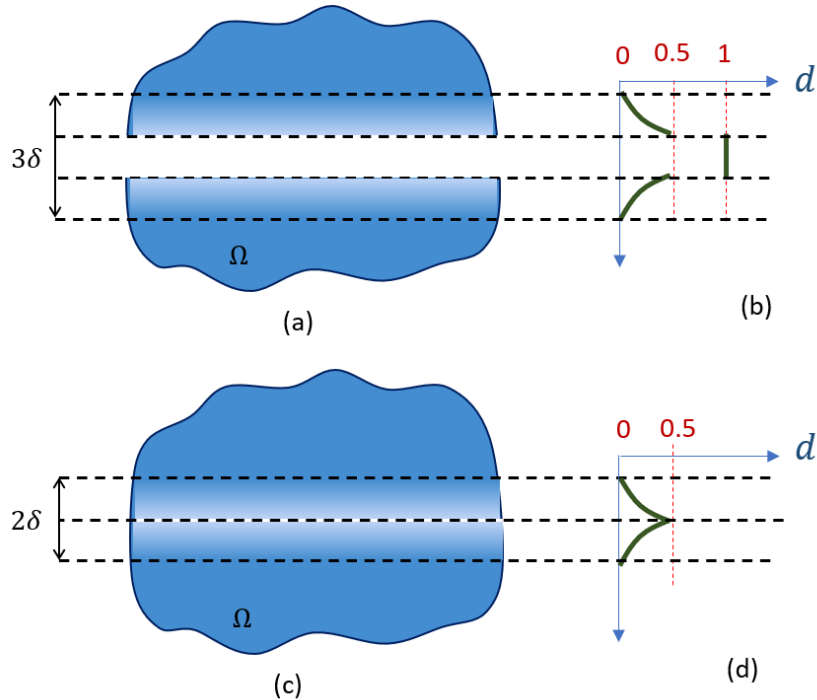

Fig. 3.   A through crack in a peridynamic body in 2D. Top is the new energy-based bond-failure model, bottom is the critical bond strain model commonly-used with the meshfree PD model.

One important observation in the model is the mass loss in the middle layer with d = 1. It is reasonable to argue that as long as $\delta$ is small relative to the geometry, the mass loss effect is negligible. In the numerical example in Section 5.2, performance of the new model is tested against the critical strain bond-breaking meshfree PD model in a dynamic fracture problem that involves multiple crack branching events.

To find a relationship between $W_c$ and the critical energy release rate $G_0$ for model calibration, we consider that all points inside the middle $\delta$-thick layer of a through crack have reached the critical strain energy $W_c$. $G_0$ is the energy released for the crack to advance a unit crack area. Let $A$ to denote the crack surface area. We enforce the equality between the strain energy required to grow the crack and create the new crack surface area $A$, and the corresponding critical energy release rate:

$$A\delta W_c = A G_0 \tag{50}$$

This results in:

$$W_c = \frac{G_0}{\delta} \tag{51}$$

We observe that this calibration is valid for brittle failure and is independent from the deformation model. Moreover, the calibrated threshold is identical for 2D and 3D.

We verify the proposed damage model in Section 5.2 by comparing the model's predictions in dynamic crack branching against published results obtained with the meshfree PD damage model, based on critical bond-strain.

### 4. The fast convolution-based method for peridynamic models in elasticity and fracture

In this section, we first discuss the convolution-based Fourier method for peridynamic models in elasticity and fracture in 3D periodic domains. Then we introduce the Embedded Constraint formulation to extend the method to problems on bounded domains of arbitrary shape and general nonlocal boundary conditions.

#### 4.1. Fourier-based method for periodic domains

Let $\boldsymbol{x} = \{x_1, x_2, x_3\}$ be the position vector in 3D and the box $\mathbb{T} = [0, L_1] \times [0, L_2] \times [0, L_3]$, be a periodic 3D domain with 0 being "identified" with $L_1, L_2,$ and $L_3$ in all three directions due to periodicity of the box. We first uniformly discretize the domain, using $N_1, N_2,$ and $N_3$ number of nodes in the Cartesian coordinate directions 1, 2, and 3 respectively:

$$\boldsymbol{x}_{nmp} = \{(n-1)\Delta x_1, (m-1)\Delta x_2, (p-1)\Delta x_3\}, \text{ where } \Delta x_1 = \frac{L_1}{N_1}; \Delta x_2 = \frac{L_2}{N_2}; \Delta x_3 = \frac{L_3}{N_3} \tag{52}$$

where $= \{1,2,\ldots,N_1\}; m = \{1,2,\ldots,N_2\}; p = \{1,2,\ldots,N_3\}$.

The total number of nodes is then $N = N_1 N_2 N_3$.

We approximate the vector-valued displacement field $\boldsymbol{u}(\boldsymbol{x}, t)$ with the following $N/2$-degree trigonometric approximation, also known as the *discrete Fourier series*. Given the discretized periodic domain [38, 39], we write:

$$u_i^N(\boldsymbol{x}_{nmp}, t) = \frac{1}{N_3 N_2 N_1} \sum_{k_3 = -\frac{N_3}{2}+1}^{\frac{N_3}{2}} \sum_{k_2 = -\frac{N_2}{2}+1}^{\frac{N_2}{2}} \sum_{k_1 = -\frac{N_1}{2}+1}^{\frac{N_1}{2}} \hat{u}_i(\boldsymbol{k}, t) e^{2\pi\zeta\left(\frac{k_1 x_1}{L_1} + \frac{k_2 x_2}{L_2} + \frac{k_3 x_3}{L_3}\right)}; i = 1,2,3 \tag{53}$$

where $\boldsymbol{k} = \{k_1, k_2, k_3\}$ is the integer vector of Fourier modes, $\zeta = \sqrt{-1}$, and:

$$\hat{u}_i(\boldsymbol{k},t) = \sum_{p=1}^{N_3}\sum_{m=1}^{N_2}\sum_{n=1}^{N_1} u_i^N(\boldsymbol{x}_{nmp},t) e^{-2\pi\zeta\left(\frac{k_1 n}{N_1}+\frac{k_2 m}{N_2}+\frac{k_3 p}{N_3}\right)}; \quad i=1,2,3 \tag{54}$$

are the *discrete Fourier coefficients* of $u_i^N$. The operation on $u_i^N$ in Eq. (54) is called the *discrete Fourier transform* (DFT), and the inverse operation on $\hat{u}_i$ in Eq.(53) is called the *inverse discrete Fourier transform* (iDFT) [38, 39]. Note that this definition of DFT maps one-to-one between $N$ values of $u_i^N$, and $N$ values of $\hat{u}_i$, meaning that the number of Fourier basis functions used in the trigonometric approximation in Eq. (40) is equal to the total number of spatial nodes.

From the computational point of view, the DFT and iDFT operations are carried out via efficient Fast Fourier Transform (FFT) algorithms, and the inverse (iFFT), which have complexity $N\log_2 N$ [16, 17].

### 4.1.1. Fourier-based discretization of peridynamics operators

To use the fast convolution-based method for PD equations, we first approximate the PD integrals via, for example, the one-point Gaussian quadrature rule (mid-point integration) on the periodic domain using the discrete nodes given by Eq. (52). Note that computing this quadrature directly, as done in the meshfree PD method [5], leads to a $O(N^2)$ complexity for each time step. Instead, for PD operators that can be expressed with the convolutional form of Eq. (14), we reorganize the quadrature summation and express it as convolution sums. To compute each convolution, we discretize the convolving functions by using their trigonometric series approximation (see Eq. (53)); we use the FFT to obtain the discrete Fourier coefficients of the convolving functions; we multiply the Fourier coefficients of these functions, and transform back the product into physical space by using the iFFT. Adding up these computed convolution terms returns the value of the PD integral.

Here we show the procedure for the PD operators discussed in Section 3:

1) The linearized bond-based elastic PD model

As described above, we first approximate the integral in Eq. (18) with Gaussian quadrature, and reorganize it to reach the convolutional structure. The convolutional forms in this section are identical to the continuous versions in Section 3, but in a discrete version. For conciseness, we denote the triple sum $\sum_{q=1}^{N_3}\sum_{s=1}^{N_2}\sum_{r=1}^{N_1}$ by $\sum_{q,s,r=1}^{N_3,N_2,N_1}$, and let $\Delta V = \Delta x_1 \Delta x_2 \Delta x_3$:

$$\begin{aligned}
L_i(\boldsymbol{x},t) \cong L_i^N(\boldsymbol{x}_{nmp},t) &= \sum_{q,s,r=1}^{N_3,N_2,N_1} C_{ij}^N(\boldsymbol{x}_{rsq}-\boldsymbol{x}_{nmp})[u_j^N(\boldsymbol{x}_{rsq},t)-u_j^N(\boldsymbol{x}_{nmp},t)]\Delta V \\
&= \sum_{q,s,r=1}^{N_3,N_2,N_1} C_{ij}^N(\boldsymbol{x}_{nmp}-\boldsymbol{x}_{rsq})u_j^N(\boldsymbol{x}_{rsq},t)\Delta V \\
&\quad - \left[\sum_{q,s,r=1}^{N_3,N_2,N_1} C_{ij}^N(\boldsymbol{x}_{nmp}-\boldsymbol{x}_{rsq})\Delta V\right] u_j^N(\boldsymbol{x}_{nmp},t) \\
&= \Delta V[C_{ij}^N \circledast_N u_j^N](\boldsymbol{x}_{nmp},t) - P_{ij}^N u_j^N(\boldsymbol{x}_{nmp},t)
\end{aligned} \tag{55}$$

where $C_{ij}^N$ is the discrete Fourier series approximation of $C_{ij}$ in Eq. (17), and $\circledast_N$ denotes the circular convolution sum (convolution sum on a periodic domain, see, e.g., [40]). $P_{ij}^N$ is a constant matrix, computable by the numerical integration of $C_{ij}^N$ using the 1-point Gaussian quadrature rule on the same discretized domain:

$$P_{ij}^N = \sum_{p,m,n=1}^{N_3,N_2,N_1} C_{ij}^N(x_{nmp}) \Delta V \tag{56}$$

Note that when we use square grids and a spherical horizon, there are nodes near the horizon edge whose volumes are only partially covered by the horizon [41, 42]. Using a fixed $\Delta V$ for all nodes in a family introduces some quadrature error as it also leaves out contributions to the integral coming from nodes just outside the horizon region but whose volumes are partially covered by it. This error can be reduced by using small grid size (large $m = \delta/\Delta x$), partial volume correction algorithms, or using decaying influence functions [42]. Partial volume correction algorithms like [41, 43], replace $\Delta V$ in Eq.(55) with a $\Delta V(|\xi|)$ which in computations can be viewed as a part of $C_{ij}^N(|\xi|)$, and therefore would preserve the convolutional structure. In the examples shown below, we use the constant $\Delta V$ algorithm, while being aware of its quadrature error.

To reduce notational complexity, we drop the superscript $N$ for discretized quantities and eliminate the argument brackets $(x_{nmp}, t)$ by replacing with the superscript $(nmp, t)$. We use $\mathbf{F}$ and $\mathbf{F}^{-1}$ to denote FFT and iFFT operations. Using FFT and iFFT, the PD integral operator can be computed as:

$$L_i^N(x_{nmp}, t) = L_i^{nmp,t} = \mathbf{F}^{-1}[\mathbf{F}(C_{ij})\mathbf{F}(u_j)]\big|^{nmp,t} \Delta V - P_{ij}u_j^{nmp,t} \tag{57}$$

Note that $j$ in the righthand side is a dummy index and denotes summations over $j = 1,2,3$ for each $i$. For demonstration, we expand Eq. (57) for $i = 1$:

$$L_1^{nmp,t} = \mathbf{F}^{-1}[\mathbf{F}(C_{11})\mathbf{F}(u_1) + \mathbf{F}(C_{12})\mathbf{F}(u_2) + \mathbf{F}(C_{13})\mathbf{F}(u_3)]\big|^{nmp,t} \Delta V - P_{11}u_1^{nmp,t} \tag{58}$$
$$- P_{12}u_2^{nmp,t} - P_{13}u_3^{nmp,t}$$

Similar relationships are obtained for $L_2^{nmp,t}$, and $L_3^{nmp,t}$ by expanding Eq. (57) for $i = 2$ and 3, and summing over the dummy index $j$.

2) The nonlinear bond-based elastic material

We follow the same discretization procedure for the PD nonlinear operator provided in Eq. (23), which yields:

$$\begin{aligned}L_i^{nmp,t} = &\mathbf{F}^{-1}[\mathbf{F}(C_{ij})\mathbf{F}(u_j)]\big|^{nmp,t} \Delta V - P_{ij}u_j^{nmp,t} - \mathbf{F}^{-1}[\mathbf{F}(a_i)\mathbf{F}(u_ju_j)]\big|^{nmp,t} \Delta V \\ &+ \mathbf{F}^{-1}[\mathbf{F}(a_i)\mathbf{F}(u_j)]\big|^{nmp,t} \Delta V u_j^{nmp,t} - \frac{1}{2}q_i^{nmp,t}u_j^{nmp,t}u_j^{nmp,t} \\ &- \mathbf{F}^{-1}[\mathbf{F}(a_j)\mathbf{F}(u_ju_i)]\big|^{nmp,t} \Delta V + u_j^{nmp,t}\mathbf{F}^{-1}[\mathbf{F}(a_j)\mathbf{F}(u_i)]\big|^{nmp,t} \Delta V \\ &+ \mathbf{F}^{-1}[\mathbf{F}(a_j)\mathbf{F}(u_j)]\big|^{nmp,t} \Delta V - q_i^{nmp,t}u_i^{nmp,t}u_j^{nmp,t} \\ &+ \frac{1}{2}\mathbf{F}^{-1}[\mathbf{F}(c)\mathbf{F}(u_ju_ju_i)]\big|^{nmp,t} \Delta V - \mathbf{F}^{-1}[\mathbf{F}(c)\mathbf{F}(u_iu_j)]\big|^{nmp,t} \Delta V u_j^{nmp,t} \\ &+ \frac{1}{2}\mathbf{F}^{-1}[\mathbf{F}(c)\mathbf{F}(u_i)]\big|^{nmp,t} \Delta V u_j^{nmp,t}u_j^{nmp,t} \\ &- \frac{1}{2}\mathbf{F}^{-1}[\mathbf{F}(c)\mathbf{F}(u_ju_j)]\big|^{nmp,t} \Delta V u_i^{nmp,t} \\ &+ \mathbf{F}^{-1}[\mathbf{F}(c)\mathbf{F}(u_j)]\big|^{nmp,t} \Delta V u_j^{nmp,t}u_i^{nmp,t} - \frac{1}{2}u_j^{nmp,t}u_j^{nmp,t}u_i^{nmp,t}p\end{aligned} \tag{59}$$

where $C_{ij}$, $a_j$, and $c$ are discrete Fourier series approximations of the functions given in Eq. (24). $P_{ij}$, $q_i$, and $p$ are the numerical integrations of $C_{ij}$, $a_j$, and $c$ respectively, similar to the quadrature in Eq. (56).

3) The linearized state-based elastic solid

According to the convolutional form derived in Eq. (38), the discrete version of $L_i(\pmb{x}, t)$ for this material is:

$$L_i^{nmp,t} = \vartheta^{nmp,t} q^{nmp,t} - \mathbf{F}^{-1}[\mathbf{F}(a_i)\mathbf{F}(\vartheta)]\Big|^{nmp,t} \Delta V + \mathbf{F}^{-1}[\mathbf{F}(C_{ij})\mathbf{F}(u_j)]\Big|^{nmp,t} \Delta V \quad (60)$$
$$- P_{ij} u_j^{nmp,t}$$

Where $C_{ij}$ and $a_j$ are discrete Fourier series approximations of the functions given in Eq.(37). $P_{ij}$ and $q_i$ are numerical integrations of $C_{ij}$ and $a_j$ respectively. The discrete version of $\vartheta$ given in Eq. (39) is:

$$\vartheta^{nmp,t} = \frac{-3}{m^{nmp,t}} \left\{ \mathbf{F}^{-1}[\mathbf{F}(a_i)\mathbf{F}(u_i)]\Big|^{nmp,t} \Delta V + q_i u_i^{nmp,t} \right\} \quad (61)$$

4) The damage model:

Here we show the discretization for the damage model introduced in Section 3.4. According to Eq. (46) one can write:

$$L_i^{nmp,t} = \lambda^{nmp,t} \left\{ \mathbf{F}^{-1}[\mathbf{F}(C_{ij})\mathbf{F}(\lambda u_j)]\Big|^{nmp,t} - \mathbf{F}^{-1}[\mathbf{F}(C_{ij})\mathbf{F}(\lambda)]\Big|^{nmp,t} u_j^{nmp,t} \right\} \Delta V \quad (62)$$

where

$$\lambda^{nmp,t} = \begin{cases} 1 & W^{nmp,t-\Delta t} \leq \dfrac{G_0}{\delta} \\ 0 & W^{nmp,t-\Delta t} > \dfrac{G_0}{\delta} \end{cases} \quad (63)$$

And $W^{nmp,t-\Delta t}$ is the nodal strain energy density at node $\pmb{x}_{nmp}$ and the previous time step $t - \Delta t$. From Eq.(47):

$$W^{nmp,t} = \frac{1}{2} \lambda^{nmp,t} \Delta V \left\{ \mathbf{F}^{-1}[\mathbf{F}(C_{ij})\mathbf{F}(u_i u_j \lambda)]\Big|^{nmp,t} - u_j^{nmp,t} \mathbf{F}^{-1}[\mathbf{F}(C_{ij})\mathbf{F}(u_i \lambda)]\Big|^{nmp,t} \right. \quad (64)$$
$$- u_i^{nmp,t} \mathbf{F}^{-1}[\mathbf{F}(C_{ij})\mathbf{F}(u_j \lambda)]\Big|^{nmp,t}$$
$$\left. + u_i^{nmp,t} u_j^{nmp,t} \mathbf{F}^{-1}[\mathbf{F}(C_{ij})\mathbf{F}(\lambda)]\Big|^{nmp,t} \right\}$$

According to Eq.(49), damage is computed by:

$$d^{nmp,t} = 1 - \frac{\lambda^{nmp,t} \mathbf{F}^{-1}[\mathbf{F}(\omega_0)\mathbf{F}(\lambda)]\Big|^{nmp,t}}{\sum_{p,m,n=1}^{N_3,N_2,N_1} \omega_0^{nmp}} \quad (65)$$

By obtaining the PD operators in discrete form, we proceed to the full discretization of the PD equation of motion for dynamic problems, and the PD equilibrium equation for static problems.

#### 4.1.2. Discretization of peridynamic equation of motion and equilibrium

Once the PD integral operator is obtained in discrete form via the FFT-based quadrature (examples shown in section 4.1.1), one can write the PD integro-differential equation of motion (see Eq. (1)) as a system of second order ODEs:

$$\rho \frac{d^2 u_i}{dt^2}\bigg|^{nmp,t} = L_i^{nmp,t} + b_i^{nmp,t} \ ; \ \text{for } i = 1,2,3 \tag{66}$$

where $b_i^{nmp,t}$ is the body force density evaluated at the node $nmp$ at time $t$. The ODE system can be solved via standard second-order ODE solvers. In the case of Velocity-Verlet time integration [13] for example:

$$v_i^{nmp,t+\frac{\Delta t}{2}} = v_i^{nmp,t} + \frac{\Delta t}{2\rho}\left(L_i^{nmp,t} + b_i^{nmp,t}\right) \tag{67}$$
$$u_i^{nmp,t+\frac{\Delta t}{2}} = u_i^{nmp,t} + \Delta t \left(v_i^{nmp,t+\frac{\Delta t}{2}}\right)$$
$$v_i^{nmp,t+\Delta t} = v_i^{nmp,t+\frac{\Delta t}{2}} + \frac{\Delta t}{2\rho}\left(L_i^{nmp,t+\Delta t} + b_i^{nmp,t+\Delta t}\right)$$

where $\boldsymbol{v}^{nmp,t} = \{v_1^{nmpt}, v_2^{nmp,t}, v_3^{nmp,t}\}$ denotes the discrete velocity field, and $\Delta t$ is the time-step.

The described Fourier-based discretization in space and Velocity-Verlet time integration can be used to solve PD equation of motion for dynamic problems in period domains. In Section 4.2 we extend the method for arbitrary domains and boundary conditions.

For static problems in periodic domains we need to solve the equilibrium equation:

$$\boldsymbol{L}(\boldsymbol{x}) + \boldsymbol{b}(\boldsymbol{x}) = 0 \tag{68}$$

which is a special case of Eq.(1), without the time dependency. Repeating the same Fourier-based discretization process given by Eqs. (52) to (66), and discarding the time variable $t$ results in:

$$L_i^{nmp} + b_i^{nmp} = 0 \ ; \ \text{for } i = 1,2,3. \tag{69}$$

Eq.(69), in general, can be a linear or nonlinear system of equations in terms of the unknown $\boldsymbol{u}^N$. We define the vector-valued *residual function* $\boldsymbol{R}(\boldsymbol{u}^N)$:

$$R_i(\boldsymbol{u}^N) = L_i^{nmp} + b_i^{nmp}; \ \text{for } \ i = 1,2,3 \tag{70}$$

One can use an iterative solver that finds $\boldsymbol{u}^N$ such that $\boldsymbol{R}(\boldsymbol{u}^N) = \boldsymbol{0}$. In this framework, when using an iterative solver that involves matrix-vector products, one should not compute such products directly in order to maintain the complexity of $O(N\log_2 N)$. The matrix-vector products in these solvers are in fact, the PD integrals operating on quantities represented by vectors. Therefore, one can use the FFT-based description of the system ($\boldsymbol{L}$ formulas given in Section 4.1.1) to compute the matrix-vector products at the cost of $(N\log_2 N)$.

The FFT-based method described in this section so far, is only applicable to the PD problems defined over periodic domains. Next, we discuss the extension of this method to bounded domains with arbitrary shapes and volume constraints using the embedded constraint approach introduced in [27] for scalar problems (diffusion), which is extended here to vector problems (elasticity).

### 4.2. Embedded constraints for enforcing boundary conditions

Before we describe the embedded constraint method we briefly discuss boundary conditions in PD nonlocal settings.

### 4.2.1. Boundary conditions in peridynamics

In the classical (local) theory, boundary conditions (BC) restrict the solution for boundary-value problems in on the boundary around the domain. These boundaries are lower dimensional manifolds compared to the domain's dimension. For example, a 3D domain has a 2D boundary, and a 2D domain has a 1D boundary.

In PD, the nonlocality requires the constraints to be defined over a "thick boundary", which has the same dimension as the domain. The constraints are defined over a chunk of the domain with the thickness of at most $\delta$ (see Fig. 4)[44].

In many applications it is usually desired/needed to apply local boundary conditions. The reason is that empirical measurements mostly provide the data on surfaces, rather than on a chunk of the domain. One approach to apply a local BC in a PD problem, is to extend the domain by a "fictitious region" to specify certain volume constraints such that the desired local BC is effectively enforced. These types of methods are usually referred to as *the fictitious nodes methods* (FNM) and are extensively discussed in [45]. There are other types of methods for applying local BC to PD problems as well (see, e.g., [46-48]).

For mechanical problems in 2D and 3D where each point has more than one degrees of freedom ($u_1$, $u_2$, $u_3$ in 3D), the constrained region is not necessarily identical for all $u_i$. For instance, where $u_1$ is constrained, $u_2$ and $u_3$ may be unknown. In this study we regard the whole domain and all constrained volume as the PD body (B). For each $u_i$, $\Omega_i$ is the domain where $u_i$ is unknown and $\Gamma_i$ is where $u_i$ is given. In the 3D case, one can write $B = \Omega_1 \cup \Gamma_1 = \Omega_2 \cup \Gamma_2 = \Omega_3 \cup \Gamma_3$. Fig. 4 shows a generic 2D PD body consisting of domains and constrained volumes.

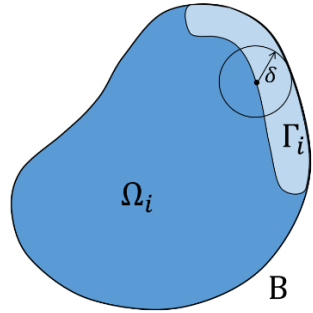

Fig. 4. A peridynamic body (B), consisting of the domain $\Omega_i$ where $u_i$ is unknown and the boundary layer/constrained volume ($\Gamma_i$) where $u_i$ is prescribed.

In this study, when we intent to use FNM to apply local BCs, the fictitious regions are regarded as $\Gamma_i$ and therefore are viewed as part of the PD body.

Dynamic problems within PD theory are often formulated as *initial-value volume-constrained PD problems* [49]:

$$\begin{cases} \rho \dfrac{\partial^2 u_i(\boldsymbol{x},t)}{\partial t^2} = L_i(\boldsymbol{x},t) + b_i(\boldsymbol{x},t) & \boldsymbol{x} \in \Omega_i, \ t > 0 \\ u_i(\boldsymbol{x},0) = u_i^0; \ v_i(\boldsymbol{x},0) = v_i^0 \ \text{(initial conditions)} & \boldsymbol{x} \in \Omega_i \\ G(u_i) = 0 \ \text{(volume constraints)} & \boldsymbol{x} \in \Gamma_i, \ t \geq 0 \end{cases} \ ; \ \text{and } i = 1,2,3 \quad (71)$$

Static problems are then defined as *volume-constrained PD problems*:

$$\begin{cases} L_i(x,t) + b_i(x,t) = 0 & x \in \Omega_i, \ t > 0 \\ G(u_i) = 0 \quad \text{(volume constraints)} & x \in \Gamma_i, \ t \geq 0 \end{cases} ; \quad \text{and } i = 1,2,3 \tag{72}$$

### 4.2.2. Embedded constraint method

In order to solve general (initial value) volume-constrained PD problems via the fast convolution-based method, the first step is to enclose the whole PD body (B) within a periodic box ($\mathbb{T}^d$) where the superscript "d" denotes the spatial dimension. Fig. 5 shows an enclosed generic body within $\mathbb{T}^2$ (compare with Fig. 4).

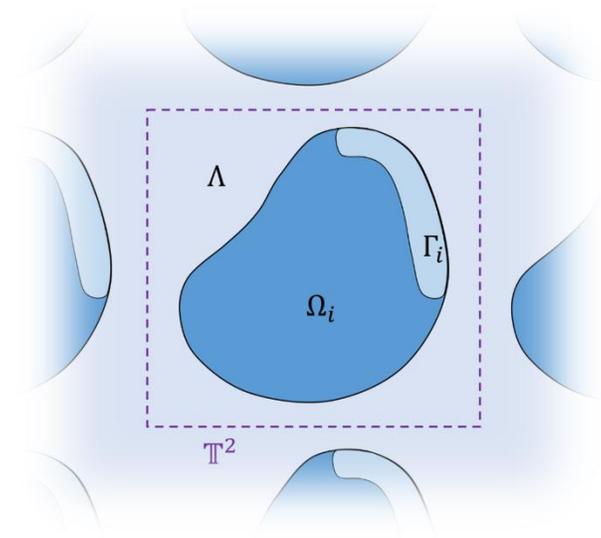

Fig. 5.　　Extension of a bounded peridynamic body to a 2D periodic box.

Note that B needs to be at least at $\delta$-distance from the edges of $\mathbb{T}^d$ to avoid any undesired nonlocal interactions between the body and the periodic extensions.

We then define the following characteristic functions:

$$\chi_B(x) = \begin{cases} 1 & x \in B \\ 0 & x \in \mathbb{T}^d \backslash B = \Lambda \end{cases} \tag{73}$$

and

$$\chi_{\Omega_i}(x) = \begin{cases} 1 & x \in \Omega_i \\ 0 & x \in \mathbb{T}^d \backslash \Omega_i = \Gamma_i \cup \Lambda \end{cases} \tag{74}$$

$\chi_B(x)$ allows us to modify the $L(x,t)$ formula for Eq. (2) and "cut" all the bonds that connect B to $\Lambda$, i.e. eliminating any interaction between the PD body and the rest of the box $\mathbb{T}^d$:

$$L_i(x,t) = \int_{\mathcal{H}_x} \chi_B(x) \chi_B(x') f_i(x,x',t) dV_{x'} = \int_{\mathcal{H}_x} \chi_B \chi_B' f_i dV_{x'} \tag{75}$$

Note that $\chi_B \chi_B'$ is always zero for bonds that have one end in B and one end in $\Lambda$. With this modification, any nonlocal interaction between the point inside and outside the body (B) is filtered by $\chi_B \chi_B' = 0$. This effectively means that the body is mechanically cutout from the rest of $\mathbb{T}^d$ and behaves as a free body.

Similar to the case for $\mu = \lambda \lambda'$ discussed in Section 3.4.2, one can show that convolutional structure of $L_i(x,t)$ is preserved after it is modified with $\chi_B \chi_B'$ in Eq. (75).

### 4.2.3. FCBM-EC for dynamic problems

In the case of PD dynamic problems consider the following initial value volume-constrained problem where $u_i(\boldsymbol{x}, t)$ are explicitly specified on $\Gamma_i$:

$$\begin{cases} \rho \dfrac{\partial^2 u_i(\boldsymbol{x},t)}{\partial t^2} = L_i(\boldsymbol{x},t) + b_i(\boldsymbol{x},t) & \boldsymbol{x} \in \Omega_i, \ t > 0 \\ u_i(\boldsymbol{x},0) = u_i^0; \ v_i(\boldsymbol{x},0) = v_i^0 \ \text{(initial conditions)} & \boldsymbol{x} \in \Omega_i \\ u_i(\boldsymbol{x},t) = g_i(\boldsymbol{x},t) \ \text{(volume constraints)} & \boldsymbol{x} \in \Gamma_i, \ t \geq 0 \end{cases} \quad ; \text{ and } i=1,2,3 \tag{76}$$

We use $\chi_{\Omega_i}(\boldsymbol{x})$ and $\chi_B(\boldsymbol{x})$ to *replace Eq.(76) which is a problem on the bounded domain B, with the following equivalent problem on the periodic domain* $\mathbb{T}^d$:

$$\begin{cases} \rho \dfrac{\partial^2 u_i}{\partial t^2} = \chi_{\Omega_i} \left( \displaystyle\int_{\mathcal{H}_x} \chi_B \chi_B' f_i dV_{x'} + b_i \right) + (1 - \chi_{\Omega_i}) \rho \dfrac{\partial^2 w_i}{\partial t^2} & \boldsymbol{x} \in \mathbb{T}^d, \ t > 0 \\ u_i(\boldsymbol{x},0) = u_i^0; \ v_i(\boldsymbol{x},0) = v_i^0 & \boldsymbol{x} \in \mathbb{T}^d \end{cases} \tag{77}$$

where $w_i(\boldsymbol{x}, t)$, is a function which is equal to volume constraints on $\Gamma_i$, and zero elsewhere:

$$w_i(\boldsymbol{x},t) = \begin{cases} g_i(\boldsymbol{x},t) & \boldsymbol{x} \in \Gamma_i \\ 0 & \boldsymbol{x} \in \mathbb{T}^d \setminus \Gamma_i \end{cases} \tag{78}$$

Since the problem in Eq. (77) is defined over a periodic domain, we can use the Fourier-based quadrature for computing $L_i^{nmp,t}$ and write:

$$\rho \left. \dfrac{d^2 u_i}{dt^2} \right|^{nmp,t} = \chi_{\Omega_i}^{nmp} \left( L_i^{nmp,t} + b_i^{nmp,t} \right) + \left(1 - \chi_{\Omega_i}^{nmp}\right) \rho \left. \dfrac{d^2 w_i}{dt^2} \right|^{nmp,t} \quad ; \text{ for } i=1,2,3 \tag{79}$$

According to the definition of $\chi_{\Omega_i}$ in Eq. (74), one can rewrite Eq. (79) as:

$$\rho \left. \dfrac{d^2 u_i}{dt^2} \right|^{nmp,t} = \begin{cases} L_i^{nmp,t} + b_i^{nmp,t}, & \boldsymbol{x}_{nmp} \in \Omega_i \\ \rho \left. \dfrac{d^2 w_i}{dt^2} \right|^{nmp,t}, & \boldsymbol{x}_{nmp} \in \mathbb{T} \setminus \Omega_i \end{cases} \quad ; \text{ for } i=1,2,3 \tag{80}$$

We use the Velocity-Verlet method for temporal integration (see Section 4.1.2) of the upper branch (on $\Omega_i$), while there is no need for temporal integration of the lower branch, since the $u_i^{nmp,t}$ values are already known (on $\Gamma_i$). One can write then:

$$u_i^{nmp,t+\Delta t} = \begin{cases} u_i^{nmp,t} + \Delta t \left[ v_i^{nmp,t} + \dfrac{\Delta t}{2\rho} \left( L_i^{nmp,t} + b_i^{nmp,t} \right) \right] & ; \ \boldsymbol{x}_{nmp} \in \Omega_i \\ w_i^{nmp,t+\Delta t} & ; \ \boldsymbol{x}_{nmp} \in \mathbb{T}^d \setminus \Omega_i \end{cases} \tag{81}$$

and

$$v_i^{nmp,t+\Delta t} \tag{82}$$
$$= \begin{cases} v_i^{nmp,t} + \dfrac{\Delta t}{2\rho} \left[ \left( L_i^{nmp,t} + b_i^{nmp,t} \right) + \left( L_i^{nmp,t+\Delta t} + b_i^{nmp,t+\Delta t} \right) \right] & ; \ \boldsymbol{x}_{nmp} \in \Omega_i \\ 0 \ \text{(value not used)} & ; \ \boldsymbol{x}_{nmp} \in \mathbb{T}^d \setminus \Omega_i \end{cases}$$

Using the characteristic function $\chi_{\Omega_i}$, Eq. (81) and (82) can be re-written as:

$$u_i^{nmp,t+\Delta t} = \chi_{\Omega_i}^{nmp}\left\{u_i^{nmp,t} + \Delta t\left[v_i^{nmp,t} + \frac{\Delta t}{2\rho}(L_i^{nmp,t} + b_i^{nmp,t})\right]\right\} \\ + (1 - \chi_{\Omega_i}^{nmp})w_i^{nmp,t+\Delta t} \qquad (83)$$

and

$$v_i^{nmp,t+\Delta t} = \chi_{\Omega_i}^{nmp}\left\{v_i^{nmp,t} + \frac{\Delta t}{2\rho}[(L_i^{nmp,t} + b_i^{nmp,t}) + (L_i^{nmp,t+\Delta t} + b_i^{nmp,t+\Delta t})]\right\} \qquad (84)$$

### 4.3. FCBM-EC for static and quasi-static problems

For static and quasi-static problems, consider the following PD volume-constrained problem where $\boldsymbol{u}(\boldsymbol{x})$ is explicitly specified on the constrained volume:

$$\begin{cases} L_i(\boldsymbol{x}) + b_i(\boldsymbol{x}) = 0 & , \boldsymbol{x} \in \Omega_i \\ u_i(\boldsymbol{x}) = g_i(\boldsymbol{x}) \text{ (volume constraints)} & , \boldsymbol{x} \in \Gamma_i \end{cases} ; \text{ and } i = 1,2,3 \qquad (85)$$

Similar to the dynamic problem, we replace this description which is defined on the bounded domain B, with the following equivalent problem on the periodic domain $\mathbb{T}^d$:

$$\chi_{\Omega_i}(\boldsymbol{x})\left[\int_{\mathcal{H}_x}\chi_B(\boldsymbol{x})\chi_B(\boldsymbol{x}')f_i(\boldsymbol{x},\boldsymbol{x}')dV_{\boldsymbol{x}'} + b_i(\boldsymbol{x})\right] + [1 - \chi_{\Omega_i}(\boldsymbol{x})][u_i(\boldsymbol{x}) - w_i(\boldsymbol{x})] = 0 \qquad (86)$$

Where

$$w_i(\boldsymbol{x}) = \begin{cases} g_i(\boldsymbol{x}) & \boldsymbol{x} \in \Gamma_i \\ 0 & \boldsymbol{x} \in \mathbb{T}\setminus\Gamma_i \end{cases} \qquad (87)$$

Discretization of Eq. (86) leads to:

$$\chi_{\Omega_i}^{nmp}(L_i^{nmp} + b_i^{nmp}) + (1 - \chi_{\Omega_i}^{nmp})(u_i^{nmp} - w_i^{nmp}) = 0; \text{ for } i = 1,2,3 \qquad (88)$$

Where $L_i^{nmp}$ is allowed to be computed using the Fourier-based quadrature due to the periodic reconstruction of the problem. Note that the second term in Eq. (88) is already satisfied by setting $u_i^N = w_i^N$. As a result we define the following residual functions:

$$R_i^{nmp}(\boldsymbol{u}^N) = \chi_{\Omega_i}^{nmp}(L_i^{nmp} + b_i^{nmp}); \text{ for } i = 1,2,3 \qquad (89)$$

One can choose an appropriate iterative solver to solve $\boldsymbol{R}(\boldsymbol{u}^N) = \boldsymbol{0}$ for $\boldsymbol{u}^N$. For the static example in Section 5.1, we use a conjugate gradient method [50] to solve for the displacement field.

### 4.4. Discussion on the accuracy and the spectral description of the method

The fast convolution based method introduced in this study for elasticity and fracture, and in [27] for diffusion problems, can be viewed as a meshfree discretization of peridynamics equations with FFT-accelerated quadrature. FCBM can also be viewed as a Fourier spectral method where the inner products with the Fourier basis functions are computed by quadrature approximation. In the Fourier spectral description of FCBM, one can assume the Fourier series expression of the convolving functions in Eq. (14) and write:

$$\int_{\mathcal{H}_x} f(\pmb{x}, \pmb{x}', t)\, dV_{\pmb{x}'} = \sum_{n=1}^{p} a_n(\pmb{x}, t)\, [b_n * c_n](\pmb{x}, t) \tag{90}$$

$$= \sum_{n=1}^{p} a_n(\pmb{x}, t)\, \mathcal{F}^{-1}[L_1 L_2 L_3 \mathcal{F}(b_n) \mathcal{F}(c_n)](\pmb{x}, t)$$

where $\mathcal{F}$ and $\mathcal{F}^{-1}$ are the exact Fourier transform and its inverse defined as:

$$\mathcal{F}(M) = \frac{1}{L_3 L_2 L_1} \int_{\mathbb{T}} M(\pmb{x}, t) e^{-2\pi \zeta \left(\frac{k_1 x_1}{L_1} + \frac{k_2 x_2}{L_2} + \frac{k_3 x_3}{L_3}\right)} dx_1\, dx_2 dx_3 \tag{91}$$

$$\mathcal{F}^{-1}(\widehat{M}) = \sum_{k_3, k_2, k_1 = -\infty}^{+\infty} \widehat{M}(\pmb{k}, t) e^{2\pi \zeta \left(\frac{k_1 x_1}{L_1} + \frac{k_2 x_2}{L_2} + \frac{k_3 x_3}{L_3}\right)}$$

In Eq. (90), if we truncate the Fourier series at $N_1, N_2, N_3$ modes in the three Cartesian directions, and approximate the $\mathcal{F}$ integrals via midpoint quadrature, using the uniform grid spacing described by Eq. (52), we recover the equations given in Section 4.1.1, having FFT-based quadrature operations.

Recall that, in general, Fourier spectral discretizations show spectral accuracy (exponential convergence rate) if the solution is smooth on $\mathbb{T}$, and if the Fourier multipliers [14] (Fourier transform of the kernels of PD operator, i.e. $\mathcal{F}(c_n)$ in Eq. (90)) are computed exactly. For example, the Fourier multipliers for the PD Laplacian operator can be expressed in terms of hypergeometric functions [14, 20]. While such analytical formulas for the Fourier multipliers lead to spectral accuracy, they depend on the kernels' forms and may not be always easy to find. Note also that evaluating these analytical relationships can be challenging (see, e.g., [20]). The quadrature approximation used in this study bounds the FCBM accuracy to that of the quadrature, but it is a general approach and can be easily used for any given kernel. Moreover, the EC method we use for incorporating boundary conditions leads to working with Fourier transforms of non-smooth functions (e.g. characteristic functions), which automatically would drop the exponential convergence rate even when analytical formulas for the Fourier multipliers are known.

The convergence rate of FCBM in terms of spatial discretization size is expected to be polynomial, bounded by quadrature's accuracy and the rate of convergence for the Fourier series approximation of transformed functions. In case of certain diffusion problems for example, the FCBM spatial rate of convergence has been shown to be quadratic [27]. We show that FCBM's spatial rate of convergence for a 3D linear elasticity example in Section 5.1 is superlinear. Rigorous error estimates and numerical analysis for FCBM is of significant interest and will be studied in the future.

Another concern with spectral methods is that aliasing errors might reduce the accuracy when a Fourier transform of a product operation is involved. We did not encounter noticeable aliasing errors in the examples shown in this study (see the convergence study in Section 5.1). However, if these errors become significant in a particular problem, one can use de-aliasing techniques such the 2/3[rd] rule [51] to remove them.

## 5. Numerical examples

In this section we solve two example problems via FCBM-EC. First, we solve a 3D example for an elastic deformation of a relatively complex geometry under static loading, using the state-based PD model. This allows us to verify our FCBM-EC static formulation, and also the convolutional form of the state-based model. Then, we solve a 2D brittle fracture problem where we verify our new damage model, as well as the FCBM-EC formulation for PD dynamic problems, with the linearized bond-based model.

### 5.1. Elastostatic deformations in a 3D body with a complex shape

In this example we compute the elastic deformation of a 3D semi-cylindrical dog bone specimen with two lateral thorough holes subjected to a static uniaxial tension. Fig. 6 shows the specimen's 3D configuration, mid-plane cross-sections and the top view, with the dimensions. The Young modulus and Poisson ratio considered are $E = 60$ GPa and $\nu = 0.4$.

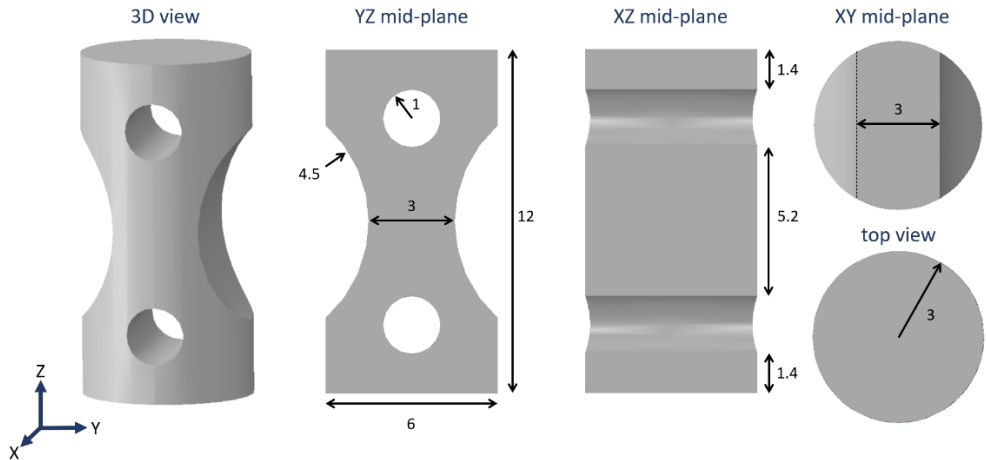

Fig. 6. The specimen's 3D configuration, YZ and XY mid-planes' cross-sectional views, and the top/bottom view. Dimensions are given in centimeters.

We select the horizon size to be $\delta = 0.3$ (sufficiently small relative to the size of the drilled holes, the smallest geometrical feature of the sample, see [52]) and extend the specimen's geometry by $\delta$ from the top and the bottom as the fictitious domain for enforcing volume constraints. These extensions are referred to as $\Gamma_3$, because $u_3$ is specified on them. Then, we further extend the PD body B (including $\Gamma_3$) in all directions by $\delta$ to form a box of size $6.6 \times 6.6 \times 14.2$ cm$^3$ as the periodic domain $\mathbb{T}^3$ (see Fig. 7).

We perform two PD simulations using two different types of nonlocal BC for each on the top and bottom boundary layers. In one simulation we set $u_3 = 0.005$ cm for all the points in $\Gamma_3$ on the top, and $u_3 = 0$ for all the points in the bottom layer. This approach is referred to as the naïve fictitious nodes method (naïve FNM) [45], and is one of the most convenient methods to determine volume constraints for approximating local BCs. In the second simulation we use the mirror-based FNM [45, 53] to define the volume constraints on $\Gamma_3$. In mirror-based FNM volume constraints on each point of $\Gamma_3$ are determined from the desired local BC value and the solution value on an interior $\Omega$ point, referred to as the "mirror" point. The description of mirror-based FNM is provided in the appendix. This approach is a better approximation of local BCs compared with the naïve FNM. However, its implementation is not as easy as the naïve approach, since the mirror-based FNM requires updating the volume constraints at each time step/iteration according to the computed solution at the previous time step/iteration. Mirror-based FNM are also more difficult to implement for curved geometries since one needs to know the normal to the surface along the boundaries. For more details please see [45].

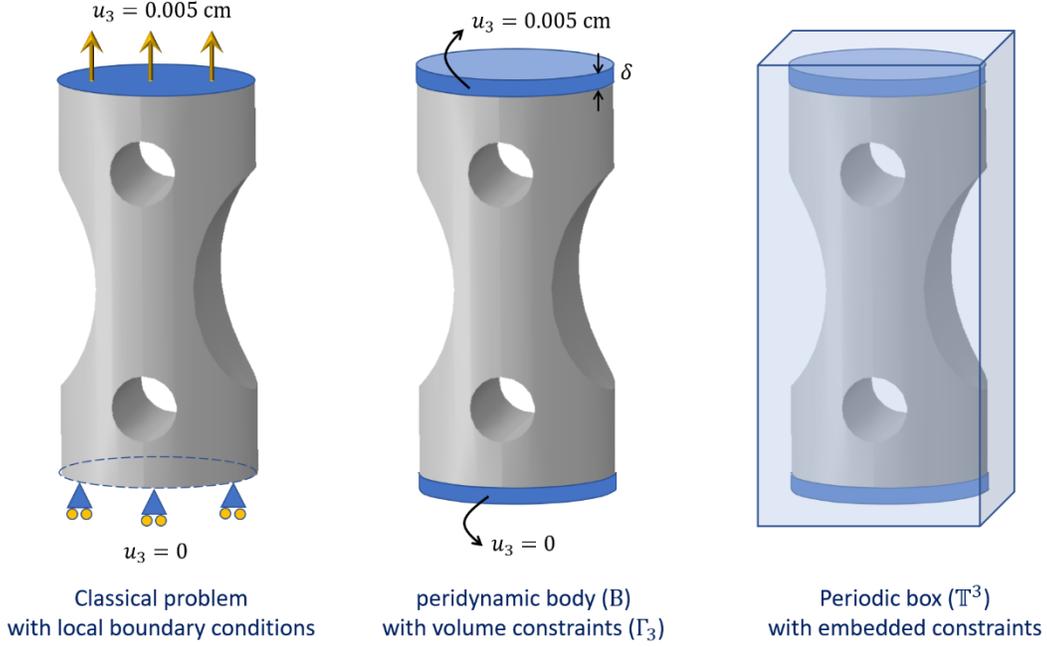

Fig. 7. Left: Specimen's configuration with local boundary conditions; Middle: PD body B with $\Gamma_3$ extensions for applying volume constraints/nonlocal BC; Right: further extension to the periodic box $\mathbb{T}^3$ for discretization with FCBM-EC .

The rest of the body's surface, which is traction-free, is treated the same in all PD simulations here, which is by eliminating nonlocal interactions between B and the rest of the box using the $\chi_B$ characteristic function (see Section 4.2.2 and Eq.(75)). This is identical to naïve FNM for traction-free conditions [45].

Considering the state-based linear elastic model in Eq. (38) and the modification in Eq. (75) to remove the nonlocal interactions between the body B with the rest of the box ($\Lambda$), the PD operator is expressed as:

$$\begin{aligned}
L_i &= \int_{\mathcal{H}_x} \chi_B \chi'_B f_i \, dV_{x'} = \int_{\mathcal{H}_x} \chi_B \chi'_B \big[a_i(x'-x)(\vartheta+\vartheta') + C_{ij}(x'-x)(u'_j - u_j)\big] \, dV_{x'} \\
&= \chi_B \bigg\{\vartheta \int_{\mathcal{H}_x} -a_i(x-x') \chi'_B dV_{x'} - \int_{\mathcal{H}_x} a_i(x-x') \chi'_B \vartheta' dV_{x'} \\
&\quad + \int_{\mathcal{H}_x} C_{ij}(x-x') \chi'_B u'_j dV_{x'} + \bigg(\int_{\mathcal{H}_x} C_{ij}(x'-x) \chi'_B dV_{x'}\bigg) u_j \bigg\} \\
&= \chi_B \big[-\vartheta(a_i * \chi_B) - (a_i * \chi_B \vartheta) + (C_{ij} * \chi_B u_j) + (C_{ij} * \chi_B) u_j \big]
\end{aligned} \tag{92}$$

Note that the disconnection between B and $\Lambda$ via $\chi_B$ must be considered in computing $\vartheta$ and $m$ as well, for consistency. Therefore, we modify Eq. (39) as follows:

$$\begin{aligned}
\vartheta &= \frac{3}{m} \int_{\mathcal{H}_x} \chi_B \chi'_B \omega(|\xi|) \xi_i \eta_i \, dV_{x'} = \frac{-3\chi_B}{m} \int_{\mathcal{H}_x} a_i(x-x') \chi'_B (u'_i - u_i) \, dV_{x'} \\
&= \frac{-3\chi_B}{m} \int_{\mathcal{H}_x} a_i(x-x') \chi'_B u'_i \, dV_{x'} - \frac{3}{m}\bigg(\int_{\mathcal{H}_x} a_i(x'-x) \chi'_B dV_{x'}\bigg) u_i \\
&= \frac{-3\chi_B}{m} \big[(a_i * \chi_B u_i) + (a_i * \chi_B) u_i\big]
\end{aligned} \tag{93}$$

where

$$m = \int_{\mathcal{H}_x} \chi_B \chi'_B \omega(|\xi|) |\xi|^2 \, dV_{x'} = \chi_B \int_{\mathcal{H}_x} \omega(x - x')|x - x'|^2 \chi'_B \, dV_{x'} = \chi_B(\omega|x|^2 * \chi_B) \tag{94}$$

Similar to $L_i$, including $\chi_B \chi'_B$ in calculation of $\vartheta$ and $m$ preserves the convolutional structure of the integrals and allows one to compute them by FFT at the cost of $O(N\log_2 N)$.

To perform the simulations, we need to choose an influence function $\omega$. Here, we use $\omega = \frac{1}{|\xi|}$ as one of the most used options [2, 54].

We discretized the domain considering $N_1 = N_2 = 2^7$ and $N_3 = 2^8$, in the three Cartesian directions. The $m$-factor (ratio of horizon to grid size) for this example happens to be the same for all directions: $\delta/\Delta x_1 = \delta/\Delta x_2 = \delta/\Delta x_3 = 5.82$. Note that the $m$-factor can be different along different directions, if grid spacing is not identical in all directions (see the example in Section 5.2).

Following the FCBM discretization described in Section 4.1.1, we can write the discretized version of the equations above as:

$$L_i^{nmp} = \chi_B^{nmp} \left\{ -\vartheta^{nmp} \mathbf{F}^{-1}[\mathbf{F}(a_i)\mathbf{F}(\chi_B)]\Big|^{nmp} - \mathbf{F}^{-1}[\mathbf{F}(a_i)\mathbf{F}(\chi_B \vartheta)]\Big|^{nmp} \right. \tag{95}$$
$$\left. + \mathbf{F}^{-1}[\mathbf{F}(C_{ij})\mathbf{F}(\chi_B u_j)]\Big|^{nmp} + \mathbf{F}^{-1}[\mathbf{F}(C_{ij})\mathbf{F}(\chi_B)]\Big|^{nmp} u_j^{nmp} \right\} \Delta V$$

where

$$\vartheta^{nmp} = -\frac{3\chi_B^{nmp}}{m^{nmp}} \left\{ \mathbf{F}^{-1}[\mathbf{F}(a_i)\mathbf{F}(\chi_B u_i)]\Big|^{nmp} + \mathbf{F}^{-1}[\mathbf{F}(a_i)\mathbf{F}(\chi_B)]\Big|^{nmp} u_i^{nmp} \right\} \Delta V \tag{96}$$

and

$$m^{nmp} = \chi_B^{nmp} \, \mathbf{F}^{-1}\left[\mathbf{F}\left(\omega(x_1^2 + x_2^2 + x_3^2)\right)\mathbf{F}(\chi_B)\right]\Big|^{nmp} \Delta V \tag{97}$$

In order to use the FCBM-EC formulation for static problems given by in Eq. (89), we construct $\chi_{\Omega_i}$ from Eq.(74), knowing that $\Omega_1 = \Omega_2 = B$, and $\Omega_3 = B \backslash \Gamma_3$. In this study, we use a conjugate gradient algorithm to solve Eq. (89).

We also solve the local version of this problem using a commercial finite element package. To this aim, we use Abaqus/Standard 6.19-1 solver with over $1 \times 10^7$ D3D4 linear tetrahedral elements and about $2 \times 10^6$ nodes, which is about the same number of degrees of freedom in $\Omega_B$ in the FCBM-EC simulation.

All the simulations in this study are performed on a Dell-Precision T7910 workstation PC, Intel(R) Xeon(R) CPU E5-2643 W v4 @3.40 GHz logical processors, and 128 GB of installed memory.

Fig. 8 shows the displacement field obtained by FCBM-EC (using the Naïve FNM boundary conditions) and by Abaqus in 3D. A rotating view of the vertical displacement field ($u_3$) obtained by FCBM is shown in Video 1.

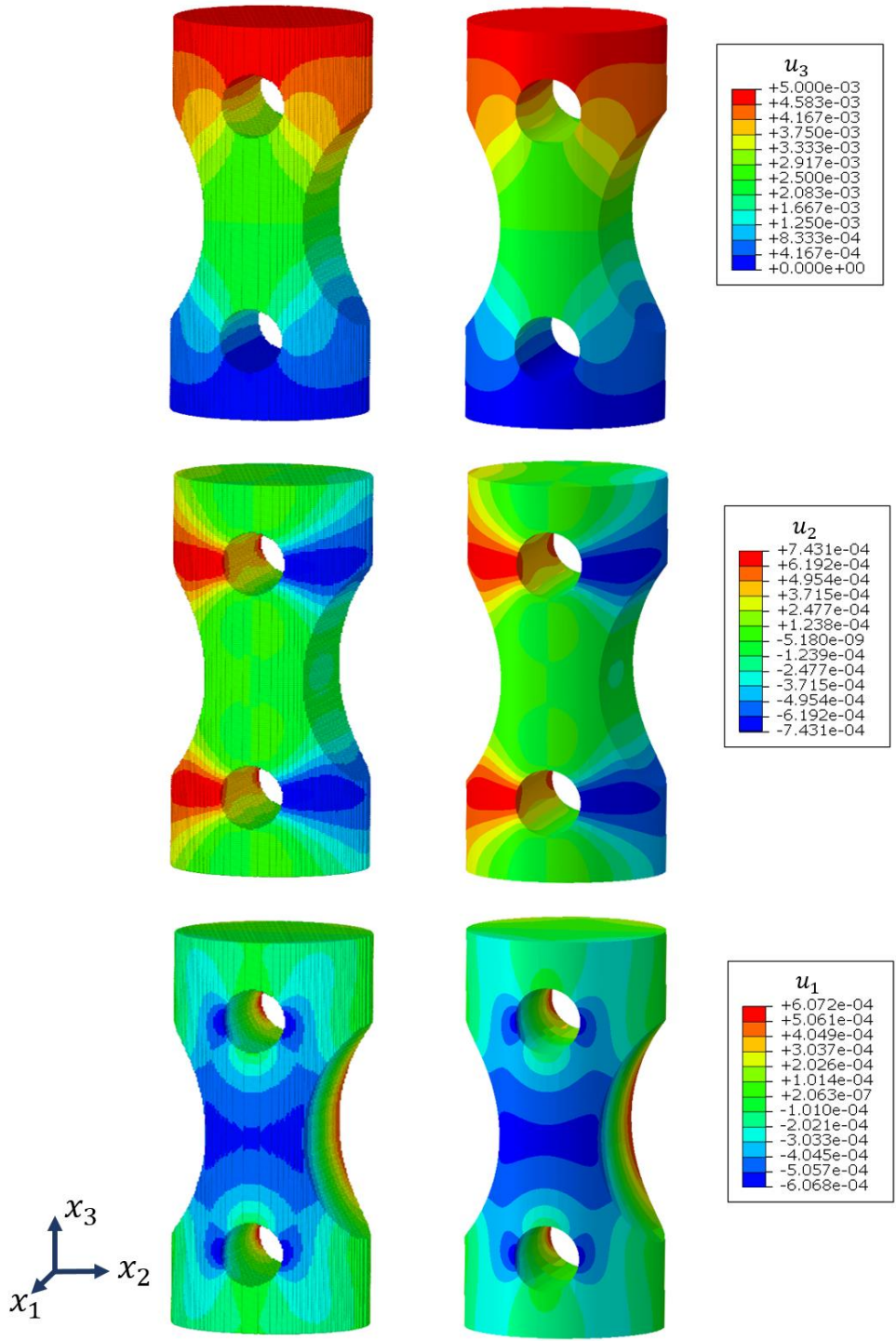

Fig. 8. Displacement field $(u_1, u_2, u_3)$ in 3D, obtained from simulations of the static tension problem using FCBM-EC discretization of peridynamics (left) and finite element solution of the classical equilibrium by Abaqus static analysis (right).

As we observe from the displacement contours, the FCBM-EC solution of the PD model and FEM solution of the classical model are visually very close. To investigate the difference quantitatively, we plot the absolute relative difference between the two solutions in the cross-sectional views in Fig. 9. The plotted differences are computed from:

$$E_i^{nmp} = \frac{\left|(u_i^{nmp})_{\text{FCBM}} - (u_i^{nmp})_{\text{Abaqus}}\right|}{\max\limits_{\text{all nodes}}(u_i^{nmp})_{\text{Abaqus}}} \tag{98}$$

where $E_i^{nmp}$ is the absolute relative difference of $u_i^{nmp}$. We intentionally use "difference" instead of "error" for referring to $E_i^{nmp}$, because the governing equations in Abaqus are different from PD (one local and one nonlocal). Therefore, the difference in the results is not only attributed to the numerical methods' error, but it also originates from the nature of the governing equations.

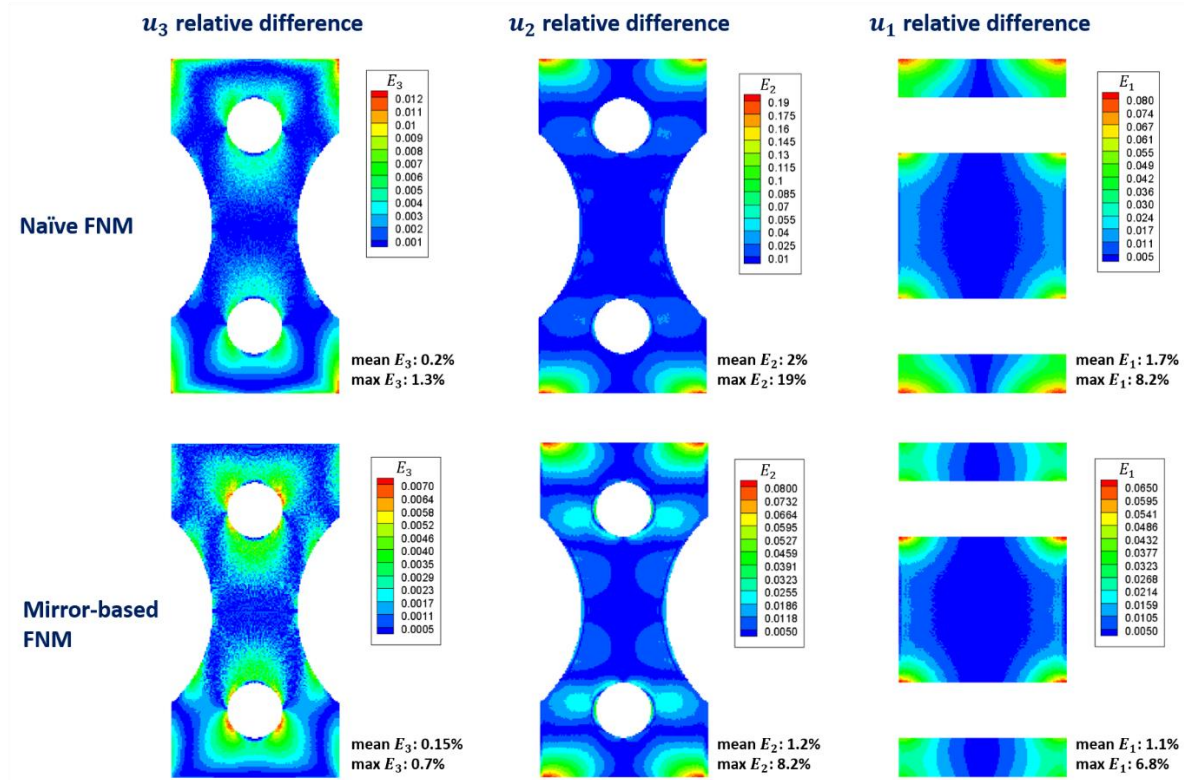

Fig. 9. Absolute relative difference plotted in cross-sectional views for displacement fields obtained via FCBM-EC in comparison with the solutions from the Abaqus static analysis (FEM solution of classical equilibrium). Top row: FCBM solver used Naïve fictitious nodes method (FNM) for boundary conditions; Bottom row: mirror-based FNM was used for the top and the bottom boundaries. Differences for $u_3$ and $u_2$ are plotted in YZ mid-plane and $u_1$ is plotted in XZ mid-plane cross-section (see Fig. 6).

The maximum difference between PD and Abaqus solutions, as expected, is observed at the corners, and this is the result of the differences between local and nonlocal BCs. Using mirror-based FNM on the top and bottom boundaries has reduced the maximum difference between the local and the nonlocal solutions associated with BCs at the corners. Mean differences for FCBM PD solutions with either FNMs, remain below 2% compared the FEM (Abaqus) solution for the classical model.

Next, we compare the computational time between the FCBM-EC and the meshfree PD method (that uses direct quadrature) for the same PD problem, using various discretization sizes. To this aim, we set $N = 2^{19}, 2^{22}, 2^{25}, 2^{28}$ for the FCBM simulation. Note that the $N$ here is the total number of nodes within the box $\mathbb{T}$. The number of nodes within the body B with these numbers are $N_B = 0.25, 2, 15.9, 127.2$ million nodes, respectively. To have a fair comparison, we use $N_B$ nodes for the meshfree PD simulations, since $N - N_B$ number of nodes locate within the gap region in FCBM-EC, which are not needed in the meshfree PD method.

Table 1 compares the computational time required to perform the simulations using FCBM-EC and the meshfree PD method. We have also taken advantage of MATLAB's built-in multi-threaded and GPU FFT functions. We performed our FCBM simulations on 1CPU, 8CPUs (using multi-threaded FFT), and on a GPU (using GPU enabled FFT) in separate tests. Multithreaded computation on more than 8 CPUs did not lead to further improvements.

Table 1. Computational time for solving the 3D example, using FCBM-EC and the meshfree PD with direct quadrature (DQ)

| *m*-factor | ~3 | ~6 | ~12 | ~23 |
|---|---|---|---|---|
| Number of nodes ($N_B$) | ~250×10³ | ~2×10⁶ | ~16×10⁶ | ~127×10⁶ |
| Meshfree PD with DQ (1 CPU) | 1 hrs | 2.8 days* | 5.9 months* | 31.1 yrs* |
| FCBM-EC (1 CPU) | 5.5 min | 53 min | 7.5 hrs | 3.9 days |
| FCBM-EC (8 CPUs) | 2 min | 16.9 min | 2.5 hrs | 1.4 days |
| FCBM-EC (GPU) | 2.7 min | 9.9 min | out of memory | out of memory |

\* time is estimated, using the time for *m*=3 and knowing that the computational time scales in O($N^2$) for the method with direct quadrature [27].

Although the extra space in the FCBM-EC computational box has almost twice the number of nodes used by the meshfree PD for the actual domain, FCBM-EC simulations are significantly faster. Due to the difference between the computational complexity of the methods, as the number of nodes increases, the efficiency gain of FCBM-EC becomes higher. The new method reduces PD computations from hours to minutes and from years to days. From Table 1, computations that were impossible to conduct using a single CPU with the meshfree PD method, are now easily achievable via FCBM-EC. Bulit-in Matlab's FFT functions allowed us to benefit from parallel computations by adding a few lines of codes in the FCBM solver. In this example, using more than 8 CPUs does not lead to much more efficiency.

Another significant advantage of FCBM-EC compared with the meshfree PD is the memory allocation. In simulations provided in Table 1, the memory allocation required by the meshfree PD only allowed for the lowest resolution ($N_B = 250 \times 10^3$). The computational time for higher spatial resolutions are estimated, knowing that they scale with O($N^2$) for the meshfree PD [27]. Most meshfree implementations initially find and store family nodes for each node to improve efficiency during the quadrature [55, 56]. Without this initialization, family search is required for each node at each time step/iteration which can potentially add a significant amount of time. However, storing family information requires allocating arrays of $N \times M$ where $M$ is the maximum number of family nodes. Knowing that for a fixed $\delta$, $M$ scales with $N$, memory allocation in the meshfree PD scales with O($N^2$). In contrast, FCBM-EC does not need to search and store family information. The integrals are computed in Fourier space, without explicitly adding individual interactions between family nodes. As a result, the variables in FCBM are arrays of size of $N \times 1$ and their storage scales as O($N$). Table 2 shows the memory allocation required by FCBM and the meshfree PD for our 3D example.

Table 2. Memory allocation required for solving the 3D example, using FCBM-EC and the meshfree PD with direct quadrature (DQ).

| $m$-factor | ~3 | ~6 | ~12 | ~23 |
|---|---|---|---|---|
| Number of nodes in the body | ~250×10³ | ~2×10⁶ | ~16×10⁶ | ~127×10⁶ |
| Meshfree PD with DQ | 5.1 GB | 327 GB* | 20.4 TB* | 1,308 TB* |
| FCBM-EC | 235 MB | 1.8 GB | 14.7 GB | 117.5 GB |

* Estimated, using the memory required for the case with $m=3$, and knowing that memory allocation scales as O($N^2$) for the method with direct quadrature [23]. The memory of the computer system here was 128 GB.

We also plot the FCBM convergence of displacement field with respect to the spatial discretization ($m$-convergence, for the horizon mentioned in the problem setup above) for this 3D elasticity problem. Since the exact analytical nonlocal solution for this 3D example is not known, we use the following error measure for the convergence study:

$$\text{relative error} = \frac{|\boldsymbol{U}_{\text{new}} - \boldsymbol{U}_{\text{old}}|_{L_2}}{|\boldsymbol{U}_{\text{new}}|_{L_2}} \tag{99}$$

where $\boldsymbol{U}$ is the nodal displacements vector: a vector containing displacements in the three directions for all nodes ($u_i^{nmp}$ for all $i,n,m,p$). The "old" and "new" subscripts denote two simulations with $\Delta x_{\text{new}} = \frac{1}{2}\Delta x_{\text{old}}$. We use the coarser grid nodal coordinates (from the "old" simulation) to compute the $L_2$ norms in Eq. (99).

In Fig. 10 we plot this relative error versus the discretization size ($\Delta x$) in log-log scale.

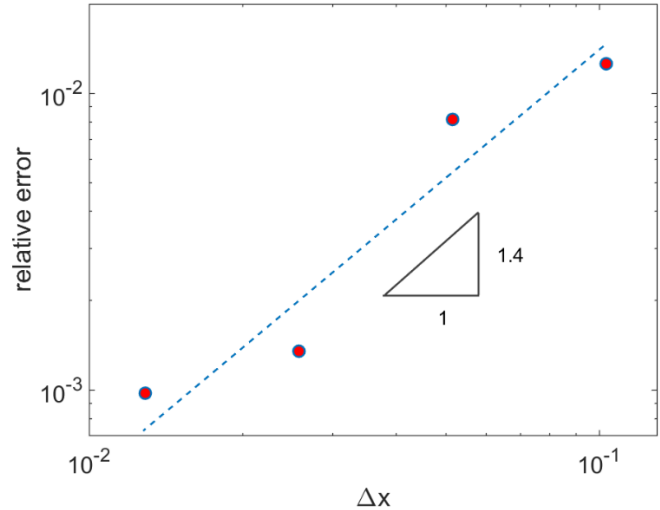

Fig. 10.   Convergence of FCBM-EC solution for displacements with respect to spatial discretization size for the 3D elasticity example. The dashed line is a linear fit of the data.

We observe a superlinear rate of convergence for this 3D elasticity example. This rate is similar to what was reported in [57], where an alternate version of this method (FCBM with volume penalization [26]) was used for a PD wave equation with discontinuous initial conditions in 2D. As discussed in Section 4.4, the accuracy of the FCBM in the present form is bounded by the quadrature and the finite Fourier series approximation errors.

In the next section, we solve a 2D fracture example and compare the results with the published meshfree PD solution to show the new method's efficiency and versatility in modeling dynamics brittle fracture problems, including those that involve multiple crack branching events.

### 5.2. Dynamic brittle fracture and crack branching in 2D

In this section we test the damage model introduced in Section 3.4 and its FCBM-EC implementation and verify it against published dynamic brittle fracture simulations obtained by the meshfree PD model that uses critical bond-strain for bond breaking [35].

The problem description is as follows: A thin plate of $10\times 4$ cm$^2$ soda-lime glass with a precrack of length 5 cm is subjected to sudden stresses that are distributed uniformly on the top and bottom boundaries and remain constant during the simulations. The sample geometry and loading conditions are shown in Fig. 11.

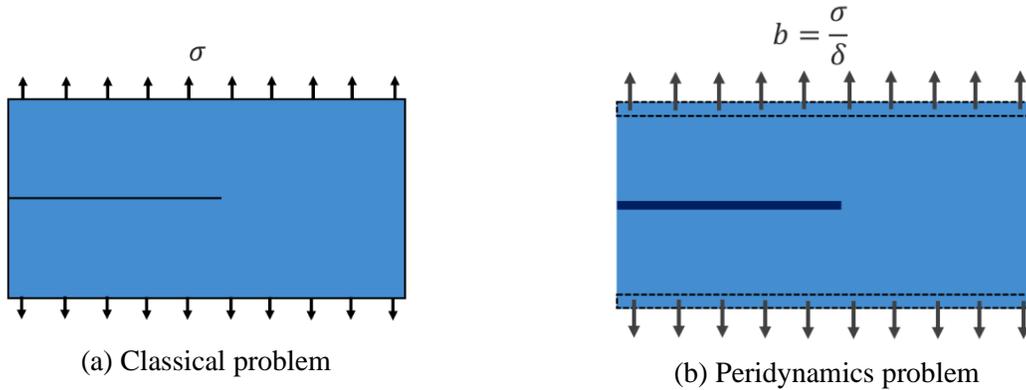

(a) Classical problem      (b) Peridynamics problem

Fig. 11. Classical problem description (a), and the corresponding peridynamics version (b).

In the peridynamic version of the problem (Fig. 11b), we impose the body force density $b = \sigma/\delta$ on the $\delta$-thick top and bottom layers of the body to enforce the desired load-controlled BCs. As noted earlier the method to impose/approximate local BCs in PD models is not unique and one may choose other approaches [45, 46, 58]. In this example, we used the body force approach, in order to have a similar loading as in [35], but our method can handle other approaches as well. Because of the through crack definition in the new damage model (see Section 3.4.3), the pre-crack in our PD model is thicker than the one used in [35] with the meshfree PD. The influence of the crack width difference in the two models is discussed after presenting the results.

The glass has the Young modulus $E = 72$ GPa, density $\rho = 2440$ Kg/m$^3$, and critical fracture energy (critical energy release rate) $G_0 = 3.8$ J/m$^2$ [35]. We conduct three simulations with three different stress values: $\sigma = 0.2, 2, 4$ MPa. These three tests are shown to result in different fracture patterns and branching behavior in [35]. We choose $\delta = 0.1$ cm as in [35].

In order to use FCBM-EC, we first extend the domain by $\delta$ in all directions to construct the periodic box (see Fig. 12).

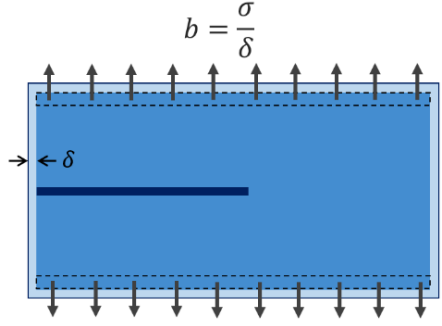

Fig. 12. Domain extension to be used in the FCBM-EC PD solver.

Note that unlike the example in Section 5.1, there are no fictitious domain in this simulations since we apply the load as a body force directly on the top and the bottom layers of the plate. Considering the damage model in Eq.(46) and the modification in Eq.(75) to disconnect the body from the rest of the gap region, the PD operator can be expressed as:

$$L_i = \int_{\mathcal{H}_x} \chi_B \chi_B' \lambda \lambda' f_i \, dV_{x'} = \chi_B \lambda \int_{\mathcal{H}_x} \chi_B' \lambda' C_{ij}(x - x')(u'_j - u_j) \, dV_{x'} \qquad (100)$$
$$= \chi_B \lambda \{[C_{ij} * (\chi_B \lambda u_j)] - u_j [C_{ij} * \chi_B \lambda]\}$$

To compute $C_{ij}$, we need to choose an influence function $\omega$ and calibrate the constant $\alpha$ to the classical elasticity constants. We use a common influence function $\omega = \frac{1}{|\xi|}$, which results in $\alpha = \frac{9E}{\pi \delta^3}$ for plane stress (Poisson ratio is restricted to $\nu = \frac{1}{3}$) [2].

We discretized the domain considering $N_1 = 2^9$ and $N_2 = 2^8$. The *m*-factor (ratio of horizon to grid size) is then $\delta/\Delta x_1 = 5.02$ and $\delta/\Delta x_2 = 6.10$ in $x_1$ and $x_2$ directions respectively. By comparing Eq. (100) with Eq.(46), one can write the discretized version by replacing $\lambda$ with $\chi_B \lambda$ in Eq. (62):

$$L_i^{nmp,t} = \chi_B^{nmp} \lambda^{nmp,t} \left\{ \mathbf{F}^{-1}[\mathbf{F}(C_{ij})\mathbf{F}(\chi_B \lambda u_j)]|^{nmp,t} \right. \qquad (101)$$
$$\left. - \mathbf{F}^{-1}[\mathbf{F}(C_{ij})\mathbf{F}(\chi_B \lambda)]|^{nmp,t} u_j^{nmp,t} \right\} \Delta V$$

To conduct the simulations, we use the FCBM-EC formulation for dynamic problems given by Eqs. (83) and (84). We choose a time step of $\Delta t = 5 \times 10^{-8}$ s which is the same as the one used by the meshfree PD simulations in [35]. In general, one can choose any $\Delta t$ that satisfies the stability condition given in [5]. At each time step $\lambda^{nmp,t+\Delta t}$ is updated using Eq.(63). This requires computing $W^{nmp,t}$. Similar to the process for obtaining the internal force density PD operator in Eq. (100), $W^{nmp,t}$ formula is found by replacing $\lambda$ with $\chi_B \lambda$ in Eq. (64):

$$W^{nmp,t} = \frac{1}{2} \chi_B^{nmp} \lambda^{nmp,t} \Delta V \left\{ \mathbf{F}^{-1}[\mathbf{F}(C_{ij})\mathbf{F}(u_i u_j \chi_B \lambda)]|^{nmp,t} \right. \qquad (102)$$
$$- u_j^{nmp,t} \mathbf{F}^{-1}[\mathbf{F}(C_{ij})\mathbf{F}(u_i \chi_B \lambda)]|^{nmp,t}$$
$$- u_i^{nmp,t} \mathbf{F}^{-1}[\mathbf{F}(C_{ij})\mathbf{F}(u_j \chi_B \lambda)]|^{nmp,t}$$
$$\left. + u_i^{nmp,t} u_j^{nmp,t} \mathbf{F}^{-1}[\mathbf{F}(C_{ij})\mathbf{F}(\chi_B \lambda)]|^{nmp,t} \right\}$$

Damage is computed from Eq. (65). Fig. 13 shows the crack patterns predicted by the new FCBM-EC PD model in comparison with the published results [35] obtained via the meshfree PD discretization of the bond-based model in Eq. (16) that uses a critical bond strain criterion, for three different loadings, at times when the cracks are about the fully split the sample.

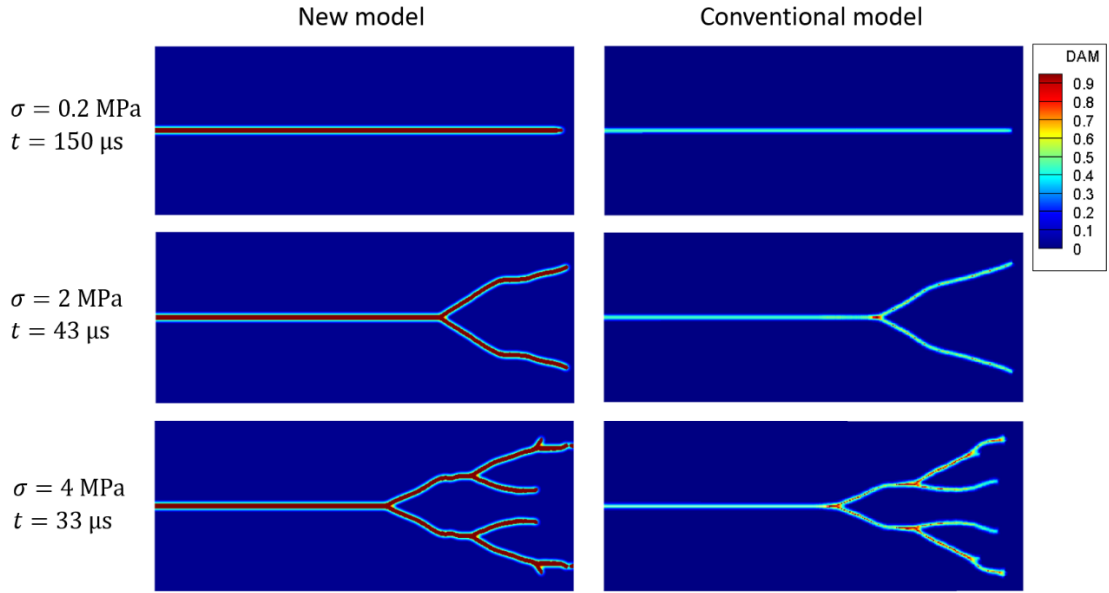

Fig. 13. Damage maps from PD simulations obtained from (left column) the new damage model (pointwise energy-based bond-breaking) solved via FCBM-EC, and from (right column) the meshfree PD solution with the critical bond-strain criterion [35].

The new damage model and solution method match, in each case, the crack patterns obtained by the meshfree PD model, including the multiple crack branching scenario. It takes the FCBM-EC model seconds to solve a problem that requires hours with the meshfree PD model.

The slight differences observed between the two simulations are expected because the two models are slightly different: the FCBM-EC discretizes the linearized bond-based PD model, while the meshfree PD solution is for the regular bond-based model; in addition, there is some difference between the damage models used. The meshfree PD simulations use the critical bond strain criterion while the FCBM solutions employ the new energy-based criterion. As discussed in Section 3.4.3, the new model shows cracks with "thickness" $3\delta$, with the middle $\delta$-thick layer having fully failed nodes ($d = 1$), whereas the meshfree PD model leads to $2\delta$-thick cracks (see Fig. 3). This difference, however, does not seems to have a significant effect on the predicted fracture patterns. From Fig. 13 it is noticed that the meshfree PD results show a slight crack thickening just before crack branching while the FCBM results show constant thickness for all crack paths. However, a closer look at the fracture kinetics at the crack tip at the time of branching in the FCBM simulation, reveals that the "physics of branching" are similar between the two methods/models. Fig. 14 shows zoom-in views of several snapshot during the branching process in the test with $\sigma = 2$ MPa for the FCBM-EC PD model.

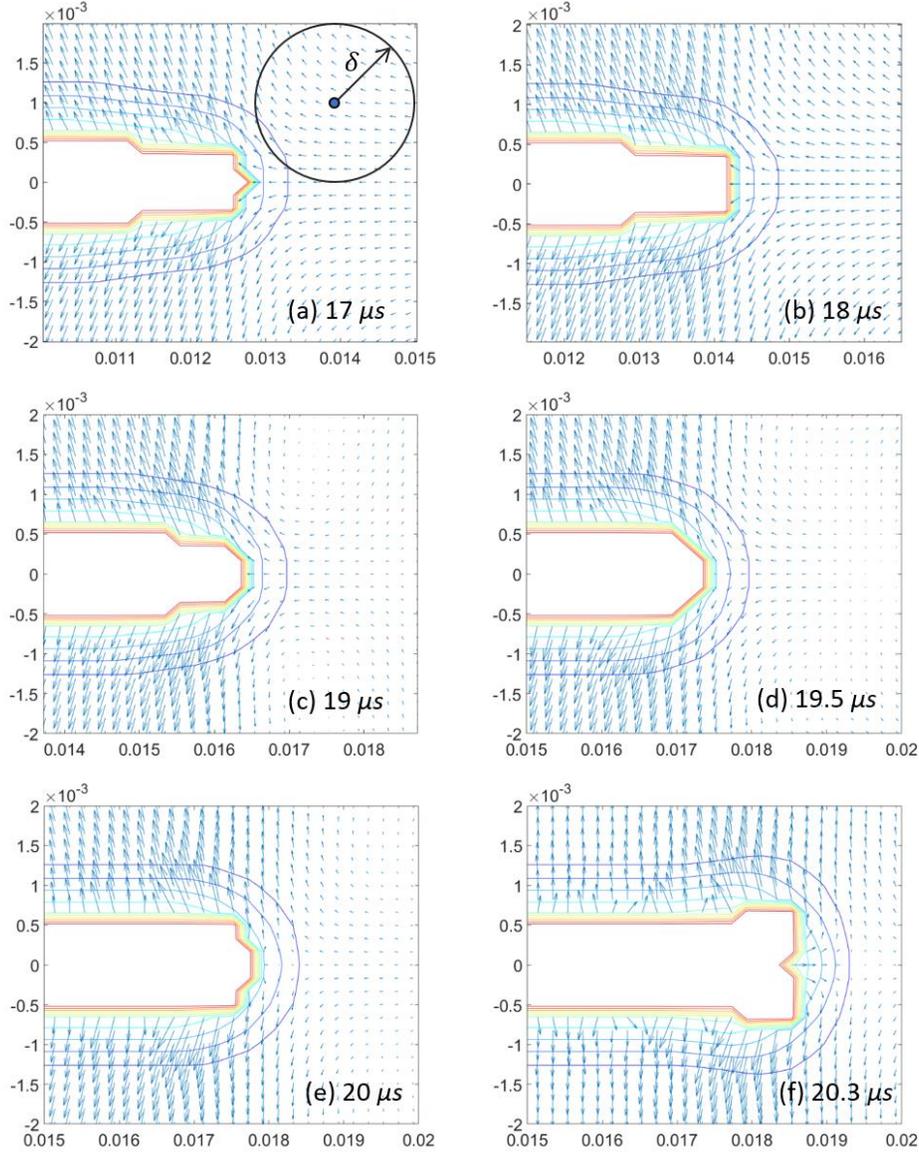

Fig. 14. Dynamics of crack branching around the crack tip with the new energy-based nodal damage model. Arrows show the velocity vector field while the contours show the damage index value. The damage contour colors are associated with the legend in Fig. 13. The PD horizon used, is plotted, to scale, in the top left picture.

We observe that the three-phase process of branching reported in [35] is observed here as well. Phase I: propagation (Fig. 14a and b); Phase II: wave pile-up and thickening of damage (Fig. 14c and d); Phase III: damage migration and branching (Fig. 14e and f). Before branching (Fig. 14a and b), while the crack grows in a straight line, damage progression happens as follows: nodes near the crack tip get damaged, which leads to the removal of all of their bonds. As the crack advances, bonds for nodes near the crack tip may still bridge the developing crack, and therefore increase the strain energy for these nodes. The strain energy density for nodes on the top and bottom crack surfaces with bridging bonds reach now the critical value and they also get removed, reaching damage index equal to one. This is why there is a widening of the crack past a horizon-length of the "process zone". This can continue as a steady state crack growth, but when the energy delivered at the crack tip is larger, then (Fig. 14c and d) the material in front of the crack tip moves more forcefully towards the advancing crack and leads to a pile-up of strain energy on the banks of the process zone, meaning that nodes on these banks reach critical strain energy not due to increased straining of their bridging bonds, but of all of their bonds: consequently, they fail sooner. This corresponds to the thickening of damage just before branching of a crack that is observed experimentally [59-61]. When this thickening reaches the tip of the process zone, (Fig. 14e and f), damage has to migrate from the center line and create the pathways for the crack to branch, since the strain energy concentration is no longer along the symmetry line, but at the corners of the thicker damage zone.

What we see, is that the FCBM PD model also captures the thickening of the crack just before branching, but this is partly obscured by the model's limitation to simulating cracks that are thinner than a thickness of $3\delta$, compared with the meshfree PD model in which the limit is $2\delta$ (see Fig. 3).

One of the important advantages of FCBM, is that one can now easily perform $m$-convergence studies (refining the grid size with a fixed horizon size) up to high values of $m$, which was previously not practical, or even impossible, with the meshfree PD method. Note that $m$-convergence studies for fracture problems are particularly important since studies have shown that a grid independent crack path is obtained only when using larger $m$ values [62] than what normally have been employed in most applications.

Here we performed an $m$-convergence study for the 2D fracture problem described above with $\sigma = 4$ MPa. We kept the horizon fixed $\delta = 0.1$ cm, and changed the number of nodes in the $x_1$ and $x_2$ directions. Six simulations were conducted with $(N_1, N_2) = (2^7, 2^6), (2^8, 2^7), (2^{10}, 2^9), (2^{10}, 2^{10})$, $(2^{11}, 2^{10})$, and $(2^{12}, 2^{11})$, leading to the $m$-factors: $(m_1, m_2) \approx (2.5, 3), (5, 6), (10, 12), (10, 24)$, $(20, 24)$, and $(40, 48)$ respectively. The case with $(m_1, m_2) \approx (10, 24)$ is chosen to test whether discretization anisotropy influences the results. Fig. 15 show the $m$-convergence results.

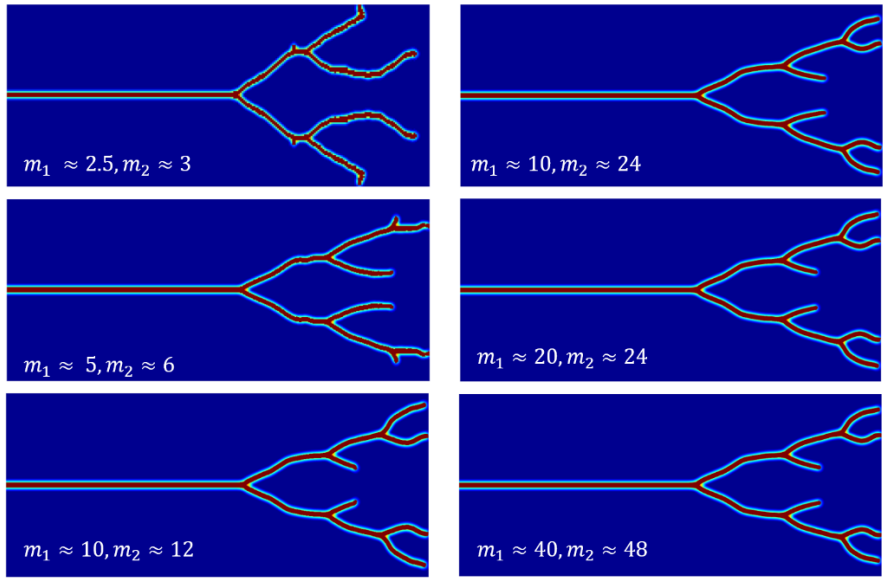

Fig. 15.   An $m$-convergence study in terms of fracture patterns for the example with suddenly applied load $\sigma = 4$ MPa.

We observe that fracture patterns converge at a value of $m \approx 10$. This is consistent with a study that concluded $m$-values higher than 7 are needed to have a correct crack path [62]. For $m$-values higher than 10, we also notice that grid spacing anisotropy does not alter the results.

In Fig. 16, we show how the simulations' run-time on a single CPU scales with the problem size.

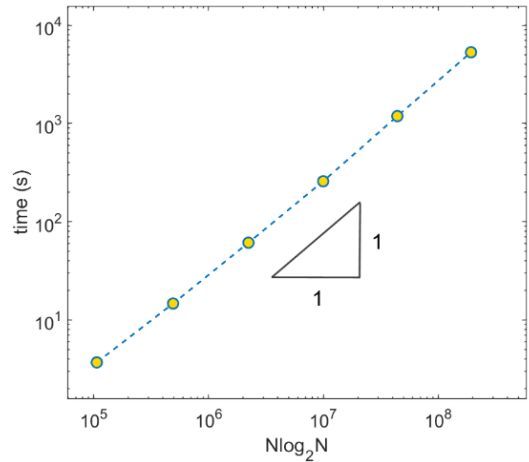

Fig. 16. Computational complexity of FCBM for this 2D fracture problem. Simulations' run-time is plotted versus $N\log_2 N$.

As expected, the simulations' time scales with $N\log_2 N$. The largest simulation here with $(m_1, m_2) \approx (40,48)$ took only 88 minutes on a single CPU, while the meshfree PD solver would have needed several months to complete the same simulation (knowing that its run-time scales with $O(N^2)$). This is an important finding: with the new model, we can now efficiently solve problems that were not accessible/feasible before: dynamic brittle fracture problems in which the fracture patterns require high *m*-values.

Evolution of displacement and velocity fields, as well as strain energy density and damage index during fracture are provided in Video 2 for the converged case with $m \approx 10$. The velocity components ($v_1$ and $v_2$) in this video, show some wave emissions created, in spurts, along the newly formed crack surfaces as the cracks grow. We investigated the discretization effect on the frequency and amplitude of these wavelets and found that these features are independent of the discretization. They are dependent then only on the damage models used, including the horizon size. Effectively, these spurts are a result of the kinetics of bond-breaking in the damage model we used here: in a dynamic process, crack propagation happens one/a few nodes at the time; i.e. one or a few nodes at the crack tip suddenly lose all of their bonds when reaching their critical strain energy density, and instantly change the force balance for their family nodes, for the next time step. This node-by-node failure process, releases kinetic energy in small packets at the crack tip, causing these small fluctuations in the velocity field.

Experimental evidence confirms that crack growth in dynamic brittle fracture happens through consecutive bursts and not in a completely smooth fashion. This is seen in fracture experiments often, as Hull states in [63], page 262: "crack growth often occurs in a series of jumps".

Similar wavelets along growing cracks, can also be observed in the videos from simulation in [35], where the meshfree PD method with the critical strain criterion was used. These wavelets are much "shaper" and better "organized" in the current FCBM PD simulations compared to the "noisier" ones from [35], and that is mainly because of the difference in the damage models: in FCBM PD, one removes a node and all of its bonds suddenly, once the node's strain energy density reaches the critical value; in the meshfree PD model, individual bonds are removed once they reach their critical value.

## 6. Conclusions

We introduced a general and fast convolution-based method (FCBM) for peridynamics. In this method, one tries to write the PD integrals in the form of convolutions. When this is possible, one then uses the one-point Gaussian quadrature over the domain and transforms the discretized convolution-based PD equations into products of discrete Fourier coefficients. Following these transformations, one computes the PD integrals with the FFT and its inverse operation at the cost of $O(N\log_2 N)$, as opposed to $O(N^2)$ in the commonly-used meshfree or FEM discretization methods of PD models. For time-dependent problems, a time-marching scheme is used. Since the integrals involved are now computed in the Fourier space, neighbor identifications and storing neighbor information is no longer needed. As a result, the memory allocation becomes of $O(N)$ in FCBM while storing neighbor information in other methods scales as $O(N^2)$. To extend the applicability of the method to arbitrary domains and boundary conditions, we extended the embedded constraint (EC) approach, previously used for PD diffusion problems, to the general setting that includes the PD equations of motion.

The method introduced here is general and can be used for any nonlocal model as long as a convolutional structure can be identified for it. In this study, we applied the procedure to the PD equations of motion and solved problems in elasticity and dynamic fracture. We showed how to obtain a convolutional structure for several material models: linearized bond- and state-based elastic materials, and a nonlinear elastic bond-based model. To exploit the FCBM efficiency for fracture problems, we introduced a new energy-based damage model that led to a convolutional structure for such problems. We tested the method on a 3D elastostatic problem over a complex geometry and a 2D dynamic brittle fracture problem with multiple crack branching events.

Comparisons between the computational efficiency of the new FCBM-EC method for PD models with that of the original meshfree discretization of PD formulations showed that problems requiring years of computations (on a single processor) with the latter method can be executed in a matter of days (on the same processor) with the former. Memory allocation was also two orders of magnitude less than what was required by the meshfree PD. One can now easily reach crack paths independent of the grid used because choosing a large number of nodes inside the PD horizon is no longer a major computational obstacle. Fast simulation of complex fracture problems with high accuracy are now possible via the FCBM-EC method for PD models, as the efficiency gains compared with the original meshfree discretization method can reach a factor of $10^3$-$10^4$ or more.

**Acknowledgments**

This work has been supported by the National Science Foundation, USA under CMMI CDS&E Award No. 1953346, and by a Nebraska System Science award from the Nebraska Research Initiative. This work was completed utilizing the Holland Computing Center of the University of Nebraska, which receives support from the Nebraska Research Initiative.